\documentclass[10pt]{amsart}
\usepackage{amsmath}
\usepackage{amscd}
\usepackage{amssymb}
\usepackage[all]{xy}
\usepackage{graphics}
\xyoption{2cell}

\newtheorem{theorem}[equation]{Theorem}
\newtheorem{proposition}[equation]{Proposition}
\newtheorem{corollary}[equation]{Corollary}
\newtheorem{lemma}[equation]{Lemma}

\theoremstyle{definition}
\newtheorem{definition}[equation]{Definition}
\newtheorem{remark}[equation]{Remark}
\newtheorem{example}[equation]{Example}
\newtheorem{definition-proposition}[equation]{Definition-Proposition}
\newtheorem{definition-corollary}[equation]{Definition-Corollary}

\numberwithin{equation}{section}

\newcommand{\arr}{\rightarrow}

\newcommand{\lgarr}{\longrightarrow}
\newcommand{\lglarr}{\longleftarrow}
\newcommand{\xarr}{\xrightarrow}

\newcommand{\cat}[1]{\operatorname{\mathsf{#1}}}
\newcommand{\iso}{\stackrel{\sim}{\rightarrow}}

\newcommand{\opn}{\operatorname}

\newcommand{\mcal}[1]{\mathcal{#1}}

\newcommand{\mrm}[1]{\mathrm{#1}}
\newcommand{\mbb}[1]{\mathbb{#1}}

\newcommand{\bsym}[1]{\boldsymbol{#1}}

\newcommand{\Hom}{\operatorname{Hom}}
\newcommand{\Ext}{\operatorname{Ext}}

\newcommand{\End}{\operatorname{End}}
\newcommand{\Proj}{\operatorname{\mathsf{Prj}}}
\newcommand{\Inj}{\operatorname{\mathsf{Inj}}}
\newcommand{\Pj}{\operatorname{\mathsf{Prj}}}
\newcommand{\pj}{\operatorname{\mathsf{prj}}}
\newcommand{\Ij}{\operatorname{\mathsf{Inj}}}
\newcommand{\Mod}{\operatorname{\mathsf{Mod}}}
\newcommand{\fmod}{\operatorname{\mathsf{mod}}}

\newcommand{\Add}{\operatorname{\mathsf{Add}}}
\newcommand{\Sum}{\operatorname{\mathsf{Add}}}
\newcommand{\ten }{\otimes}

\newcommand{\lten}{{\otimes}^{\boldsymbol{L}}}

\newcommand{\thick}{\operatorname{\mathsf{thick}}}
\newcommand{\tri}{\operatorname{\mathsf{tri}}}
\newcommand{\Loc}{\operatorname{\mathsf{Loc}}}
\newcommand{\Morph}{\operatorname{\mathsf{Mor}}}

\newcommand{\z}{\mathbb{Z}}
\newcommand{\nt}{\operatorname{\mathsf{noeth}}}

\newcommand{\bo}{\operatorname{b}}
\newcommand{\sm}{\operatorname{sm}}
\newcommand{\se}{\operatorname{se}}
\newcommand{\co}{\operatorname{c}}

\renewcommand{\AA}{\mcal{A}}
\newcommand{\BB}{\mcal{B}}
\newcommand{\CC}{\mcal{C}}

\newcommand{\MM}{\mcal{M}}
\newcommand{\RR}{\mcal{R}}
\renewcommand{\SS}{\mcal{S}}
\newcommand{\TT}{\mcal{T}}
\newcommand{\UU}{\mcal{U}}
\newcommand{\VV}{\mcal{V}}
\newcommand{\PP}{\mcal{P}}
\newcommand{\CCC}{\cat{C}}
\newcommand{\KKK}{\cat{K}}
\newcommand{\DDD}{\cat{D}}
\def\KK#1{\KKK_N(\mcal{#1})}
\def\DD#1{\DDD_N(\mcal{#1})}

\usepackage{color}

\title[Derived categories of $N$-complexes]
{Derived categories of $N$-complexes}
\author{Osamu Iyama, Kiriko Kato and Jun-ichi Miyachi}
\date{\today}
\address{O. Iyama: Graduate School of Mathematics, Nagoya University Chikusa-ku, Nagoya, 464-8602
Japan}
\email{iyama@math.nagoya-u.ac.jp}
\address{K. Kato:  Graduate School of Science, Osaka Prefecture University,
1-1 Gakuen-cho, Nakaku, Sakai, Osaka 599-8531, JAPAN}
\email{kiriko@mi.s.osakafu-u.ac.jp}
\address{J. Miyachi: Department of Mathematics, Tokyo Gakugei
University, Koganei-shi, Tokyo, 184-8501, Japan}
\email{miyachi@u-gakugei.ac.jp}

\begin{document}
\maketitle

\begin{abstract}
We study the homotopy category $\KKK_{N}(\BB)$ of $N$-complexes of an additive category $\BB$
and the derived category $\DDD_{N}(\AA)$ of an abelian category $\AA$.
First we show that both $\KKK_N(\mcal{B})$ and $\DDD_N(\mcal{A})$ have natural structures of triangulated categories. 
Then we establish a theory of projective (resp., injective) resolutions and derived functors.
Finally, under some conditions of an abelian category $\AA$,  we show that $\DDD_{N}(\AA)$ is  
triangle equivalent to 
the ordinary derived category $\DDD(\Morph_{N-2}(\AA))$ where 
$\Morph_{N-2}(\AA)$ is the category of sequential $N-2$ morphisms of $\AA$.
\end{abstract}


\setcounter{section}{-1}
\section{Introduction}\label{intro}

The notion of $N$-complexes, that is, graded objects with $N$-differentials $d$
($d^{N}=0$), was introduced by Mayer \cite{Ma} in his study of simplicial complexes.
Recently Kapranov and Dubois-Violette gave abstract framework of homological theory of $N$-complexes \cite{Ka,D}. 
Since then the $N$-complexes attracted many authors, for example 
\cite{Be,BDW,CSW,DK,Gi,GH,HK,Ka,Mi1,Mi2}.
The aim of this paper is to give a solid foundation of homological
algebra of $N$-complexes by generalizing classical theory of derived
categories due to Grothendieck-Verdier.
In particular we study homological algebra of $N$-complexes of an abelian category $\AA$ based on the modern point
of view of Frobenius categories (see \cite{Ha} for the definition) and their corresponding algebraic
triangulated categories.

In section \ref{NcpxI}, we study the category $\CCC_N(\BB)$ of $N$-complexes 
over an additive category $\mathcal B$ and the homotopy category $\KKK_N(\BB)$. 
Precisely speaking, we introduce an exact structure on $\CCC_N(\BB)$ to prove the following results.

\begin{theorem}[Theorems \ref{NcpxFrob} and \ref{NcpxTricat}]
\begin{enumerate}
\item
The category $\CCC_N(\BB)$  has a structure of a Frobenius category.
\item
The category $\KKK_N(\BB)$  has a structure of a triangulated category.
\end{enumerate}
\end{theorem}

We give an explicit description of the suspension functor $\Sigma$ and
triangles in $\KKK_N(\BB)$. 
Unlike the classical case $N=2$, the suspension functor $\Sigma$ does not coincide
with the shift functor $\Theta$. However we have the following connection between
$\Sigma$ and $\Theta$ in $\KKK_N(\BB)$.

\begin{theorem}[Theorem \ref{suspshft}]
There is a functorial isomorphism $\Sigma^2\simeq\Theta^N$ on $\KKK_N(\BB)$.
\end{theorem}

In Section \ref{DNcpx}, we introduce the derived category $\DDD_N(\AA)$ of $N$-complexes
for an abelian category $\AA$. We generalize the theory of projective resolutions
of complexes initiated by Verdier \cite{Ve} and extended to unbounded complexes by Spaltenstein and B\"ockstedt-Neeman  \cite{Sp,BN}.
Our main result is the following, where $\Proj \AA$ (resp., $\Inj \AA$) is the subcategory
of projective (resp., injective) objects in $\AA$ and 
$\KKK^{\mrm{a}}_{N}(\AA)$ (resp., $\KKK^{\mrm{p}}_{N}(\AA)$, $\KKK^{\mrm{i}}_{N}(\AA)$) is the homotopy category of $N$-acyclic
(resp., $\KKK$-projective, $\KKK$-injective) $N$-complexes (see Definitions \ref{acyclic}, \ref{dfn:spcpx}). 
We denote by $\KKK_{N}^{-}(\cat{Prj}\AA)$ (resp., $\KKK_{N}^{-,\bo}(\cat{Prj}\AA)$, $\KKK_{N}^{-,\mrm{a}}(\cat{Prj}\AA)$)
the subcategory of $\KKK_{N} (\cat{Prj}\AA)$ consisting of $N$-complexes bounded above (resp., bounded above with bounded homologies, bounded above and $N$-acyclic).
For other unexplained notations, we refer to the paragraph before Theorem \ref{cor:subeqv}. 

\begin{theorem}[Theorems \ref{cor:subeqv} and \ref{cor:subeqv03}]
The following hold for $\natural=$nothing$,\bo$.
\begin{enumerate}
\item Assume that $\AA$ has enough projectives.
\begin{enumerate}
\item $(\KKK_{N}^{-,\natural}(\cat{Prj}\AA), \KKK_{N}^{-,\mrm{a}}(\AA))$ is a stable t-structure
in $\KKK_{N}^{-,\natural}(\AA)$ and we have triangle equivalences
$\KKK_{N}^{-}(\cat{Prj}\AA) \simeq\DDD_{N}^{-}(\AA)$ and
$\KKK_{N}^{-,\bo}(\cat{Prj}\AA) \simeq\DDD_{N}^{\bo}(\AA)$.
\item If $\AA$ is an $Ab4$-category, then 
$(\KKK_{N}^{\mrm{p}}(\AA), \KKK_{N}^{\mrm{a}}(\AA))$ is a stable t-structure
in $\KKK_{N}(\AA)$ and we have a triangle equivalence
$\KKK_{N}^{\mrm{p}}(\AA) \simeq\DDD_{N}(\AA)$.
\end{enumerate}
\item Assume that $\AA$ has enough injectives.
\begin{enumerate}
\item $(\KKK_{N}^{+,\mrm{a}}(\AA),\KKK_{N}^{+,\natural}(\cat{Inj}\AA))$ is a stable t-structure
in $\KKK_{N}^{+,\natural}(\AA)$ and we have triangle equivalences
$\KKK_{N}^{+}(\cat{Inj}\AA) \simeq\DDD_{N}^{+}(\AA)$
and $\KKK_{N}^{+,\bo}(\cat{Inj}\AA) \simeq\DDD_{N}^{\bo}(\AA)$.
\item If $\AA$ is an $Ab4^*$-category, then 
$(\KKK_{N}^{\mrm{a}}(\AA),\KKK_{N}^{\mrm{i}}(\AA))$ is a stable t-structure
in $\KKK_{N}(\AA)$ and we have a triangle equivalence
$\KKK_{N}^{\mrm{i}}(\AA) \simeq\DDD_{N}(\AA)$.
\end{enumerate}
\end{enumerate}
\end{theorem}

Moreover, 
we generalize a result of Krause \cite{Kr2} characterizing the compact
objects in classical homotopy categories. We deal with a \emph{locally noetherian
Grothendieck category}, that is, a Grothendieck category with a set of generators
of noetherian objects. We give the following result, where $\CC^{\co}$ denotes the
subcategory of compact objects in an additive category $\CC$.

\begin{theorem}[Theorem \ref{Kijcp2}]
Let $\AA$ be a locally noetherian Grothendieck category
with the subcategory $\cat{noeth} \AA$ of noetherian objects in $\AA$.
\begin{enumerate}
\item
$\KKK_N(\Inj\AA)$ is compactly generated.
\item
The canonical functor $\KKK_N(\Inj \AA )\to \DDD_N(\AA)$ induces an equivalence
between $\KKK_{N}(\Ij \AA)^{\co}$
and $\DDD_{N}^{\bo}(\cat{noeth} \AA)$.
\end{enumerate}
\end{theorem}

We generalize the classical existence theorem of derived functors
to our setting by showing that any triangle functor $\KKK_N(\AA)\to \KKK_{N'}(\AA')$ has a
left/right derived functor $\DDD_N(\AA)\to \DDD_{N'}(\AA')$ (see Definition \ref{derf}) under certain
mild conditions on $\AA$. Our result is the following.

\begin{theorem}[Theorem \ref{thm:exderfun02}]
Let $\AA$, $\mcal{A'}$ be abelian categories, 
$F : \KKK_{N}(\AA) \arr \KKK_{N'}(\mcal{A'})$ a triangle functor. 
Then the following hold.
\begin{enumerate}
\item If $\AA$ is an $Ab4$-category  with enough projectives, then
the left derived functor $\bsym{L}F:\DDD_{N}(\AA) \arr \DDD_{N'}(\mcal{A'})$
exists.
\item If $\AA$ is an $Ab4^*$-category  with enough injectives, then
the right derived functor $\bsym{R}F:\DDD_{N}(\AA) \arr \DDD_{N'}(\mcal{A'})$
exists.
\end{enumerate}
\end{theorem}

In section \ref{TrieqDN}, we give our main result in this paper. We show that
the derived category $\DDD_{N}(\AA)$ is triangle equivalent to
the ordinary derived category $\DDD(\Morph_{N-2}(\AA))$ of 
$\Morph_{N-2}(\AA)$, where $\Morph_{N-2}(\AA)$ is the category 
of sequences of $N-2$ morphisms of $\AA$ (see Definition \ref{smcat}).

\begin{theorem}[Theorems \ref{DND} and \ref{DND2}]\label{realize1}
Let $\AA$ be an $Ab3$-category with a small full subcategory of compact projective generators.
Then we have a triangle equivalence for $\natural=$nothing$,+,-,\bo$.
\[
\DDD_{N}^{\natural}(\AA) \simeq\DDD^{\natural}(\Morph_{N-2}(\AA)).
\]
\end{theorem}

As applications, we have the following triangle equivalences. 
Here $\BB$ is an additive category, $\Morph^{\sm}_{N-2}(\BB)$ is the category 
of sequences of $N-2$ split monomorphisms of $\BB$ (see Definition \ref{smcat}) 
and $T_{N-1}(R)$ is the upper triangular matrix ring of size $N-1$ over a ring $R$. 
For a full subcategory $\mcal{C}$ of an additive category $\BB$ with arbitrary coproducts,
$\Sum_{\BB}{\mcal{C}}$ is the category of direct summands of coproducts of objects of $\mcal{C}$ in $\BB$.
For a ring $R$, $\fmod R$ (resp., $\pj R)$ is the category of finitely presented (resp., finitely generated projective) $R$-modules.

\begin{corollary}[Corollary \ref{KNhtp}, Proposition \ref{DNMod}]\label{realize2}
\begin{enumerate}
\item Let $\BB$ be an additive category with arbitrary coproducts.
If the subcategory $\BB^{\co}$ of compact objects of $\BB$ is skeletally small and satisfies $\BB=\Sum(\BB^{\co})$,
then we have triangle equivalences $\KKK^{-}_{N}(\BB)\simeq \KKK^{-}(\Morph^{\sm}_{N-2}(\BB))$ and
$\KKK^{\bo}_{N}(\BB)\simeq\KKK^{\bo}(\Morph^{\sm}_{N-2}(\BB))$.
\item For a ring $R$, we have a triangle equivalence $\KKK^{\natural}_{N}(\pj R) \simeq \KKK^{\natural}(\pj \opn{T}_{N-1}(R))$ for $\natural=-, \bo, (-,\bo)$.
For a right coherent ring $R$, we have a triangle equivalence $\DDD^{\natural}_{N}(\fmod R)\simeq \DDD^{\natural}(\fmod \opn{T}_{N-1}(R))$ for $\natural=$nothing$,-, \bo$.
\end{enumerate}
\end{corollary}

In \cite{IKM3}, we will study more precise relations between the homotopy categories. 

\section{Preliminaries}

In this section, we collect preliminary results on additive and triangulated categories. We will omit proofs of elementary facts.

\begin{lemma}\label{prop:pullpush3}
In an abelian category, consider a pull-back (resp., push-out) diagram 
\[\xymatrix@R=1em{
X\ \ar[d]_{f} \ar[r]^{g} & X' \ar[d]^{f'} \\
Y\ \ar[r]^{g'} & Y'}\]
and morphisms $\left(\begin{smallmatrix}g' & f' \end{smallmatrix}\right): X'\oplus Y \to Y'$, $\left(\begin{smallmatrix}g \\ f \end{smallmatrix}\right): X \to X'\oplus Y$.
Then the following hold.
\begin{enumerate}
\item  If $f'$ (resp., $f$) is epic (resp., monic), then 
the above diagram is also push-out (resp., pull-back), 
and $f$ (resp., $f'$) is also epic (resp., monic).
\item  The induced morphism $\opn{Ker}f \arr \opn{Ker}f'$ is an isomorphism
(resp., an epimorphism).
\item  The induced morphism $\opn{Cok}f \arr \opn{Cok}f'$ is a monomorphism
(resp., an isomorphism).
\item We have an exact sequence
$
0 \to \opn{Cok}f \to \opn{Cok}f' \to \opn{Cok}\left(\begin{smallmatrix}g' & f' \end{smallmatrix}\right) \to 0\ (\mbox{resp.,}\\ \ 0 \to \opn{Ker}\left(\begin{smallmatrix} g \\ f \end{smallmatrix}\right) \to \opn{Ker} f \to \opn{Ker}f' \to 0
$.
\end{enumerate}
\end{lemma}

A commutative square is called \emph{exact} if it is pullback and push-out \cite{Po}.

\begin{lemma}\label{pullback2}
In an abelian category, consider two pull-back squares (X) and (Y)
\[\xymatrix@R1em{
\ar @{} [dr] |{(X)} A \ar[d]_{a} \ar[r]& \ar @{} [dr] |{(Y)} B \ar[d]^{b} \ar[r] &C \ar[d]^{\co} \\
D \ar[r] & E \ar[r] & F.
}\]
Then the square (X+Y) is exact if and only if the squares (X) and (Y) are exact.
\end{lemma}

\begin{lemma}\label{split square} 
In an abelian category, consider an exact square with a split epimorphism $d$.
\[\xymatrix@R1.5em{
A\oplus B\ \ar[d]_{\iota=\left(\begin{smallmatrix}\iota _1 & \iota _2 \end{smallmatrix}\right)} \ar[rr]^{\left(\begin{smallmatrix}0&1\end{smallmatrix}\right)} &&
B\ar[d]^{\left(\begin{smallmatrix}1\\ 0\end{smallmatrix}\right)} \\
D\ \ar[rr]^{d=\left(\begin{smallmatrix}d_1 \\ d_2 \end{smallmatrix}\right)} && B\oplus C }\]
Then there exists an isomorphism $a:A\oplus B\oplus C\to D$ such that
$\iota =a\left(\begin{smallmatrix}1&0\\ 0&1\\ 0&0\end{smallmatrix}\right)$
and $da=\left(\begin{smallmatrix}0&1&0\\ 0&0&1\end{smallmatrix}\right)$.
\end{lemma}

\begin{proof}
Since $d$ is a split epimorphism, there exists $\iota_3:C\to D$ such that $d_1\iota_3=0$ and $d_2\iota_3=1$.
Then $a=(\iota_1\ \iota_2\ \iota_3)$ satisfies the desired conditions.
\end{proof}

For a triangulated category $\TT$ and a full subcategory $\CC$ of $\TT$, we denote by
$\tri\CC=\tri_{\TT}\CC$ the smallest triangulated subcategory of $\TT$ containing $\CC$,
and by $\thick\CC=\thick_{\TT}\CC$ the smallest triangulated subcategory of $\TT$ containing $\CC$ and closed under direct summands, 
and by $\Loc\CC=\Loc _{\TT}\CC$ the smallest triangulated subcategory of $\TT$ containing $\CC$ and closed under coproducts.
\begin{definition}[Triangle Functor] \label{dfn:dfun}
Let $\TT$ and $\TT'$ be triangulated categories 
with suspensions $\Sigma_{\TT}$ and $\Sigma_{\TT'}$ respectively. 
A {\it triangle functor}  is a pair $(F, \alpha)$, where
$F:\TT \arr \TT'$ is an additive functor 
and $\alpha : F\Sigma_{\TT} \iso \Sigma_{\TT'}F$ is a functorial isomorphism 
such that $(FX, FY, FZ, F(u), F(v), {\alpha}_XF(w))$ is a  
triangle in $\TT'$ whenever
$(X, Y, Z, u, v, w)$ is a triangle in $\TT$.
If a triangle functor $F$ is an equivalence, then we say that
$\TT$ is {\it triangle equivalent} to $\TT'$.

Let $(F,\alpha), (G,\beta) :\TT \arr \TT'$ be triangle functors. A {\it functorial morphism of triangle functors} is a functorial morphism $\phi : F \arr G$ satisfying
$(\Sigma_{\TT'}\phi)\alpha = \beta\phi \Sigma_{\TT}$.
\end{definition}

Let $\TT$ be a triangulated category and $\UU$, $\VV$ be full subcategories. 
The category of extensions $\UU* \VV$ is the full subcategory of $\TT$ 
consisting of objects $X$ such that there exists a triangle 
$U \to X \to V \to \Sigma U$ with $U \in \UU$ and $V \in \VV$. 

Note that $(\UU*\VV )*\mcal{W} = \UU*(\VV*\mcal{W})$ holds by octahedral axiom. 

\begin{definition}[\cite{Miy}]\label{st-tprime}
Let $\TT$ be a triangulated category.
A pair $(\UU, \VV)$ of full triangulated subcategories of $\TT$ is called a {\it stable t-structure} (also known as
\emph{semiorthogonal decomposition}, \emph{torsion pair}, \emph{Bousfield localization}) in $\TT$ provided that
\[\opn{Hom}_{\TT}(\UU, \VV) = 0\ \mbox{ and }\ \TT=\UU*\VV.\]
\end{definition}
In this case, the canonical quotient $\TT\to\TT/\UU$ (resp., $\TT\to\TT/\VV$) has a right (resp., left) adjoint,
and we have a triangle equivalence $\TT / \UU \simeq \VV$ (resp., $\TT / \VV \simeq \UU$).

\begin{lemma}\label{fully faithful}\cite{JK}
Let $\TT$ be a triangulated category and $\UU$, $\VV$ be full triangulated subcategories. Then the following conditions are equivalent.
\begin{enumerate}
\item $\VV*\UU\subset\UU*\VV$.
\item $\UU*\VV$ is a triangulated subcategory of $\TT$.
\item Any morphism $f:U\to V$ with $U\in\UU$ and $V\in\VV$ factors through an object in $\UU\cap\VV$.
\end{enumerate}
In this case, $(\UU/ (\UU \cap \VV) , \VV/ (\UU \cap \VV) )$ is a stable t-structure in
$(\UU*\VV )/ (\UU \cap \VV)$. Hence 
we have triangle equivalences 
$\UU/ (\UU \cap \VV ) \simeq (\UU*\VV )/ \VV$ and $\VV/ (\UU \cap \VV ) \simeq (\UU*\VV ) / \UU$.
Thus the canonical functors 
$\UU/(\UU\cap\VV)\to\TT/\VV$ and $\VV/(\UU\cap\VV)\to\TT/\UU$ are fully faithful. 
\end{lemma}

\section{Homotopy category of $N$-complexes}\label{NcpxI}

In this section, we study the homotopy category of $N$-complexes. 
We fix a positive integer $N\ge2$.
Throughout this section $\BB$ is an additive category.
An \emph{$N$-complex} $X=(X^i, d_{X}^i)$ is a diagram
\[\cdots \xrightarrow{d_X^{i-1}}X^i\xrightarrow{d_X^i}X^{i+1}\xrightarrow{d_X^{i+1}}\cdots\]
with $X^i\in\BB$ and $d_X^i\in\Hom_{\BB}(X^i,X^{i+1})$ satisfying
\[d_{X}^{i+N-1}\cdots d_{X}^{i+1}d_{X}^{i}=0\]
for any $i\in\z$.
We often denote the $r$-th power of $d_X$ by
\[d_X^{\{r\}}=d_{X}^{i+r}\cdots d_{X}^{i+1}d_{X}^{i}\]
 without mentioning grades, where $d_X^{\{0\}}=1$.
A \emph{morphism} $f:X \to Y$ between $N$-complexes is a commutative diagram
\[\xymatrix@R1.5em{
\cdots\ar[r]^{d_X^{i-1}}&X^i\ar[r]^{d_X^i}\ar[d]^{f^i}&X^{i+1}\ar[r]^{d_X^{i+1}}\ar[d]^{f^{i+1}}&\cdots\\
\cdots\ar[r]^{d_Y^{i-1}}&Y^i\ar[r]^{d_Y^i}&Y^{i+1}\ar[r]^{d_Y^{i+1}}&\cdots}\]
with $f^i\in\Hom_{\BB}(X^i,Y^i)$ for any $i\in\z$.
We denote by $\CCC_{N}(\BB)$ the category of $N$-complexes.

We call an $N$-complex $X$ \emph{bounded above} (resp., \emph{bounded below}) if 
$X^i=0$ for all $i \gg 0$ (resp., $i \ll 0$), 
and \emph{bounded} if $X$ is both bounded above and bounded below. We denote by
$\CCC_{N}^{-}(\BB)$ (resp., $\CCC_{N}^{+}(\BB)$,
$\CCC_{N}^{\bo}(\BB)$) the full subcategory of bounded above (resp., bounded below,
bounded) $N$-complexes.

Our approach to the category $\CCC_{N}(\BB)$ of $N$-complexes is based on the theory of exact categories
\cite{Qu} (see \cite{Ke0} for modern account).
Let $\SS_N(\BB)$ be the collection of sequences $0\to X\xrightarrow{f}Y\xrightarrow{g}Z\to0$ of 
morphisms in $\CCC_{N}(\BB)$ such that $0\to X^i\xrightarrow{f^i}Y^i\xrightarrow{g^i}Z^i\to0$ is split exact
in $\BB$  for any integer $i$.
Then we have the following basic observation.

\begin{theorem}\label{NcpxFrob}
The category $(\CCC_N(\BB), \SS_N(\BB))$ of $N$-complexes is a Frobenius category.
\end{theorem}

For an object $M$ of $\BB$ and integers  $s$ and $1\leq r \leq N$, let 
\[
\mu^{s}_{r}(M): \cdots \to 0 \to M^{s-r+1} \xarr{d^{s-r+1}} \cdots \xarr{d^{s-2}} M^{s-1} \xarr{d^{s-1}} M^{s} \to 0 \to \cdots
\] 
be an $N$-complex given by $M^{s-i}=M$ ($0 \leq i \leq r-1$) and $d^{s-i}=1_{M}$ ($0 < i \leq r-1$).  
One can easily check the functorial isomorphisms 
{\small
\begin{equation}\label{adjoint} 
\opn{Hom}_{\CCC_N(\mathcal{B})}(X, \mu_N^{s}(M))  \simeq \opn{Hom}_{\mathcal{B}}(X^{s}, M) \mbox{ and }
\opn{Hom}_{\CCC_N(\mathcal{B})}(\mu_N^{s}(M), X) \simeq \opn{Hom}_{\mathcal{B}}(M,X^{s-N+1})
\end{equation} }
where $f \in\opn{Hom}_{\mathcal{B}}(X^{s}, M)$ and $g \in\opn{Hom}_{\mathcal{B}}(M,X^{s-N+1})$ are mapped to 
$\rho^s_f$ and $\lambda^s_g$ respectively by the following commutative diagrams. 
{\small\[\xymatrix@R1.5em{
\mu_N^{s}(M)\ar[d]^{\rho_{f}^{s}}: &\cdots\ar[r]&0\ar[r]\ar[d] & M \ar[r]^{1}\ar[d]^{f} & \cdots \ar[r]^{1} & M\ar[d]^{d^{\{N-1\}}f}\ar[r] &0\ar[r]\ar[d]&\cdots \\
X\ar[d]^{\lambda_{g}^{s}}: & \cdots \ar[r]^d & X^{s-N} \ar[r]^d\ar[d] &X^{s-N+1} \ar[r]^d\ar[d]^{gd^{\{N-1\}}} & \cdots \ar[r]^d & X^s \ar[r]^d\ar[d]^g & X^{s+1} \ar[r]^d\ar[d] & \cdots \\
\mu_N^{s}(M): &\cdots\ar[r]&0\ar[r] & M \ar[r]^{1} & \cdots \ar[r]^{1} & M\ar[r]&0\ar[r]&\cdots \\
}\]}
\begin{lemma}\label{NcpxFrob03}
The object $\mu_{N}^{s}(M)$ is projective-injective in $(\CCC_N(\BB), \SS_N(\BB))$
for any object $M \in \BB$ and  any integer $s$.
\end{lemma}

\begin{proof}
For any exact sequence $0 \to X \to Y \to Z \to 0$ in
$\mathcal{S}_N(\mathcal{B})$, the isomorphism \eqref{adjoint} gives a commutative diagram of exact sequences
{\small\[\xymatrix@R1.5em{
0 \ar[r] & \opn{Hom}_{\CCC_N(\mathcal{B})}(Z, \mu_N^{s}(M))  \ar[r]\ar[d]^{\wr} &  \opn{Hom}_{\CCC_N(\mathcal{B})}(Y, \mu_N^{s}(M))  \ar[r] \ar[d]^{\wr} &
 \opn{Hom}_{\CCC_N(\mathcal{B})}(X, \mu_N^{s}(M)) \ar[d]^{\wr} \\
 0 \ar[r] & \opn{Hom}_{\mathcal{B}}(Z^{s}, M)  \ar[r] &  \opn{Hom}_{\mathcal{B}}(Y^{s}, M)  \ar[r] &
 \opn{Hom}_{\mathcal{B}}(X^{s}, M)  \ar[r] & 0,}
\]}
where the lower sequence is exact since $0\to X^s\to Y^s\to Z^s\to0$ is split exact.
This means that $\mu_{N}^{s}(M)$ is injective. Dually one can show that $\mu_{N}^{s}(M)$ is projective. 
\end{proof}

Let $X\in\CCC_N(\BB)$ be given. 
We have morphisms $\rho_{1_{X^{n-N+1}}}^n:\mu_N^n(X^{n-N+1}) \to X$ and $\lambda_{1_{X^{n}}}^n:X\to\mu_N^n(X^{n})$, 
using \eqref{adjoint}. 
Set $\rho_X=(\rho_{1_{X^{n-N+1}}}^n)_n:\bigoplus_{n \in \mathbf{Z}}\mu_N^n(X^{n-N+1}) \to X$ and
$\lambda_X=(\lambda_{1_{X^{n}}}^n)_n:X\to \bigoplus_{n\in \mathbf{Z}}\mu_N^n(X^{n})$.
Then we have the following exact sequences in $\mathcal{S}_N(\mathcal{B})$.
{\small\begin{equation}\label{hull and cover}
0 \to \opn{Ker}\rho_X\xarr{\epsilon_X} \bigoplus_{n \in \mathbf{Z}}\mu_N^n(X^{n-N+1}) \xarr{\rho_X} X \to 0,\quad
0 \to X \xarr{\lambda_X}\bigoplus_{n\in \mathbf{Z}}\mu_N^n(X^{n}) \xarr{\eta_X} \opn{Cok}\lambda_X  \to 0.
\end{equation}}

\begin{proof}[of Theorem \ref{NcpxFrob}]
The exact sequences \eqref{hull and cover} with Lemma \ref{NcpxFrob03} show that $(\CCC_N(\BB), \SS_N(\BB))$ has enough projectives and enough injectives.
Let $X$ be an arbitrary projective (resp., injective) object. Then, on the first (resp., second) sequence of \eqref{hull and cover}, $X$ is a direct summand of the middle term. 
By Lemma \ref{NcpxFrob03}, $X$ is injective (resp., projective).
\end{proof}

The \emph{stable category} $\underline{\mcal{F}}$ of a Frobenius category $(\mcal{F}, \SS)$
has the same objects as $\mcal{F}$ and the homomorphism set between $X,Y\in\underline{\mcal{F}}$ is given by
\[\Hom _{\underline{\mcal{F}}}(X, Y )=\Hom_{\mcal{F}}(X,Y)/ \mcal{I}(X,Y)\]
where $\mcal{I}(X,Y)$ is the subgroup of $\Hom _{\mcal{F}}(X, Y )$
consisting of morphisms which factor through some projective-injective object of $(\mcal{F}, \SS)$.
By \cite{Ha}, $\underline{\mcal{F}}$ has a structure of a triangulated category, which is nowadays called an \emph{algebraic triangulated category}.

Now we shall describe the stable category of our Frobenius category $(\CCC_N(\BB), \SS_N(\BB))$ more explicitly. 
Indeed, as in the classical case, it coincides with 
the homotopy category of $N$-complexes.
Recall that a {morphism} $f:X \to Y$ of $N$-complexes
is called \emph{null-homotopic} if there exists $s^i\in\Hom_{\mathcal{B}}(X^i,Y^{i-N+1})$
such that
\begin{equation}\label{eq:null homotopic}
f^i=\sum_{j=1}^{N-1}d_Y^{i-1}\cdots d_Y^{i-N+j}s^{i+j-1}d_X^{i+j-2}\cdots d_X^i
\end{equation}
for any $i\in\z$.
For morphisms $f,g:X \to Y$ in $\CCC_N(\BB)$, we denote $f \sim g$ if
$f-g$ is null-homotopic.
We denote by $\KKK_{N}(\BB)$ the \emph{homotopy category}, that is, the category consisting of
$N$-complexes such that the homomorphism set between $X,Y \in \KKK_{N}(\BB)$ is given by
\[\Hom_{\KKK_{N}(\BB)}(X,Y)=\Hom_{\CCC_{N}(\BB)}(X,Y)/\sim.\]

\begin{theorem}\label{NcpxTricat}
The stable category of the Frobenius category $(\CCC_N(\BB),\SS_N(\BB))$ is the homotopy category $\KKK_N(\BB)$ of $\BB$.
In particular, $\KKK_N(\BB)$ is an algebraic triangulated category.
\end{theorem}

\begin{proof}
It suffices to show that a morphism $f:X \to Y$ is null-homotopic if and only if
$f$ factors through the morphism $\lambda_X:X \to \bigoplus_{n\in \mathbf{Z}}\mu_N^n(X^{n})$ given in \eqref{hull and cover}.
This can be easily checked by \eqref{adjoint}.
\end{proof}

Now we define functors $\Sigma,\Sigma^{-1}:\CCC_{N}(\BB)\to\CCC_{N}(\BB)$ by
\[\Sigma^{-1}X=\opn{Ker}\rho_X\ \mbox{ and }\ \Sigma X=\opn{Cok} \lambda_X\]
in the exact sequences \eqref{hull and cover}. Then $\Sigma$ and $\Sigma ^{-1}$ induce the suspension functor and its quasi-inverse of the triangulated category
$\KKK_N(\BB)$.

On the other hand, we define the \emph{shift functor} $\Theta:\CCC_{N}(\BB) \to \CCC_{N}(\BB)$ by
\[
\Theta(X)^{i} = X^{i+1}\ \mbox{ and }\ d_{\Theta(X)}^{i} = d_X^{i+1}
\]
for $X=(X^{i}, d_X^{i}) \in \CCC_{N}(\BB)$.
This induces the shift functor $\Theta:\KKK_N(\BB)\to\KKK_N(\BB)$ 
which is a triangle functor. 
Unlike classical case, $\Sigma$ does not coincide with $\Theta$. However we have the following observation.

\begin{theorem}\label{suspshft}
There is a functorial isomorphism 
$\Sigma^2 \simeq \Theta^N$ on $\KKK_{N}(\BB)$.
\end{theorem}

To prove this, we give a more explicit description of $\Sigma$ and $\Sigma^{-1}$.
Let $X= (X^i,d^i)$ be an object of $\CCC_{N}(\BB)$. In
\eqref{hull and cover}, 
the first sequence is given by
{\footnotesize 
\[ (\Sigma ^{-1} X)^m = \bigoplus _{i=m-N+1}^{m-1} X^i, \quad 
d_{\Sigma ^{-1} X}^m = \left( \begin{array}{c|ccccc}
-d&1&0&\cdots&0&0\\
-d^{\{2\}}&0&1&\cdots&0&0\\
\vdots&\vdots&\vdots&\ddots&\vdots&\vdots\\
-d^{\{N-3\}}&0&0&\cdots&1&0\\
-d^{\{N-2\}}&0&0&\cdots&0&1\\
\hline
-d^{\{N-1\}}&0&0&\cdots&0&0\\
\end{array} \right) \] 
\[ (\epsilon _X )^m = \left( \begin{array}{ccccc}
1&0&0&\cdots&0\\
-d&1&0&\cdots&0\\
0&-d&1&\cdots&0\\
\vdots&\vdots&\ddots&\ddots&\vdots\\
0&0&\cdots&-d&1\\
0&0&\cdots&0&-d \\
\end{array} \right) ~\mbox{and}~
(\rho _X )^m = \left ( \begin{array}{ccccc}
d^{\{N-1\}} &d^{\{N-2\}} & \cdots & d& 1 
\end{array} \right).\]}
while the second sequence by
{\footnotesize
\[ (\Sigma X)^m = \bigoplus _{i=m+1}^{m+N-1} X^i, \quad 
d_{\Sigma X}^m = \left( \begin{array}{c|ccccc}
0&1&0&\cdots&0&0\\
0&0&1&\cdots&0&0\\
\vdots&\vdots&\vdots&\ddots&\vdots&\vdots\\
0&0&0&\cdots&1&0\\
0&0&0&\cdots&0&1\\
\hline
-d^{\{N-1\}}&-d^{\{N-2\}}&-d^{\{N-3\}}&\cdots&-d^{\{2\}}&-d\\
\end{array} \right), \] 

\[ (\lambda _X )^m  = \left ( \begin{array}{c}
1\\ d\\ \vdots\\ d^{\{N-2\}}\\ d^{\{N-1\}}\\
\end{array} \right)  ~\mbox{and}~
(\eta_X )^m = \left( \begin{array}{cccccc}
-d&1&0&\cdots&0&0\\
0&-d&1&\cdots&0&0\\
0&0&-d&\cdots&0&0\\
\vdots&\vdots&\vdots&\ddots&\ddots&\vdots\\
0&0&0&\ldots&-d&1 \\
\end{array} \right) . \]
}

\begin{proof}[of Theorem \ref{suspshft}]
We shall construct a functorial isomorphism $\Sigma \to \Theta ^{N} \Sigma ^{-1}$. 
Given an object $X=(X^i,d^i)\in \CCC_N(\BB)$, we have
$(\Sigma  X )^m = \bigoplus _{i =m+1}^{m+N-1} X^i = (\Sigma ^{-1} X) ^{m+N}$
for each $m$ by \eqref{hull and cover}. 
Let $\phi _X ^m : (\Sigma X )^m \to (\Sigma ^{-1} X) ^{m+N}$ be a morphism given as 
{\footnotesize \[ \phi _X ^m = \left( \begin{array}{ccccc}
1&0&0&\cdots&0\\
d&1&0&\cdots&0\\
d^{\{2\}}&d&1&\cdots&0\\
\vdots&\vdots&\ddots&\ddots&\vdots\\
d^{\{N-2\}}&d^{\{N-3\}}&\cdots&d&1\\
\end{array} \right).\] }
Then it is easy to check that $\phi _X$ makes the following diagram commutative
\[ \xymatrix@R1.5em{ 
(\Sigma X)^m \ar[d]^{\phi _X  ^m} \ar[rr]^{d_{\Sigma  X}^m } && (\Sigma X)^{m+1} \ar[d]^{\phi _X ^{m+1}}\\
(\Sigma ^{-1}  X)^{m+N} \ar[rr]^{d_{\Sigma ^{-1}  X}^{m+N} } && (\Sigma ^{-1}  X)^{m+N+1}.\\
} \]
Thus $\phi_X:\Sigma X\to\Theta^N\Sigma^{-1}X$ is an isomorphism in $\CCC_N(\BB)$.
\par\noindent
Next let $f$ be a morphism from $X$ to $Y$ in $\CCC_N(\BB)$. 
It is routine to show $(\Theta^N \Sigma ^{-1} f )   \phi _X = \phi _Y \Sigma f $ holds.
Thus $\phi$ gives a functorial isomorphism $\Sigma \simeq \Theta ^N\Sigma^{-1}$. 
\end{proof}

We denote by $\KKK_{N}^{-}(\BB)$ (resp., $\KKK_{N}^{+}(\BB)$, $\KKK_{N}^{\bo}(\BB)$) 
the full subcategory of $\KKK_{N}(\BB)$ corresponding to $\CCC_{N}^{-}(\BB)$
(resp., $\CCC_{N}^{+}(\BB)$, $\CCC_{N}^{\bo}(\BB)$).
Then they are full triangulated subcategories of $\KKK_N(\BB)$ by the above descriptions of $\Sigma$ and $\Sigma^{-1}$.

\begin{definition}[Hard truncations]\label{hard truncations}
For an $N$-complex $X = (X^{i}, d^{i})$, set
\[\begin{aligned}
{\tau}_{\leq n}X & : \cdots \arr X^{n-2} \arr X^{n-1}\arr X^{n} \arr 0 \arr \cdots ,\\
{\tau }_{\geq n}X & : \cdots \arr 0 \arr X^{n}\arr X^{n+1}\arr X^{n+2}\arr \cdots .
\end{aligned}\]
Then we have a triangle $\tau_{\ge n}X\to X\to\tau_{\le n-1}X\to\Sigma(\tau_{\ge n}X)$ in $\KKK_N(\BB)$.
\end{definition}

Later we will use the following observation.

\begin{lemma}\label{Sigma of mu}
We have the following.
\begin{enumerate}
\item For any $C\in\BB$, $i,s\in\z$ and $0<r<N$, we have
$\Sigma^{2i+k}\mu_{r}^{s}(C)\simeq\left\{\begin{array}{ll}
\mu_{r}^{-iN+s}(C )&(k=0)\\
\mu_{N-r}^{-iN+s-r}(C)&(k=1).
\end{array}\right.$
\item $\KKK_N^{\bo}(\BB)=\tri\{\mu^s_1(C)\mid C\in\BB,\ 0<s<N\}$. 
\end{enumerate}
\end{lemma}

\begin{proof}
(1) For each $C \in \BB$ and $r,i \in \z$ with $1 \leq r \leq N-1$, we have a term-wise split exact sequence $0 \to \mu_{r}^{-iN+s}(C) \to \mu_{N}^{-iN+s}(C) \to \mu_{N-r}^{-iN+s-r}(C) \to 0$ in $\CCC(\BB)$.
Since $\mu_{N}^{-iN+s}(C)$ is a projective-injective object in $\CCC_{N}(\BB)$, we have the desired isomorphisms in $\KKK_{N}(\BB)$.
\par\noindent
(2) Using triangles in Definition \ref{hard truncations}, we can show $\KKK_N^{\bo}(\BB)=\tri\{\mu^s_1(C)\mid C\in\BB,\ s\in\z\}$ by an induction on the number of non-zero terms. Moreover, we can replace the condition $s\in\z$ by $0\le s<N$ since $\Sigma^2\simeq\Theta^N$ holds by Theorem \ref{suspshft}. We can further replace it by $0<s<N$ since $\mu^0_1(C)=\Sigma\mu^{N-1}_{N-1}(C)$ belongs to $\tri\{\mu^s_1(C)\mid C\in\BB,\ 0<s<N\}$.
\end{proof}

We end this section with an explicit description of the mapping cone.
{For} a morphism  $f:Y=(Y^i, e^i)\to X=(X^i, d^i)$ in $\CCC_N(\BB)$, the mapping cone $\opn{C}(f)$ is given by the  diagram
\[ \xymatrix@R1.5em{ 
0 \ar[r] &Y \ar[r]^{ \lambda _Y} \ar[d]^{f} & I(Y) \ar[d]^{\psi_f} 
\ar[r]^{\eta_Y} &\Sigma Y  \ar[r]  \ar@{=}[d]  & 0\\
0 \ar[r] &X \ar[r]^{g}  & \opn{C}(f)  \ar[r]^{h} &{\Sigma Y} \ar[r]   & 0, \\ } \] 
{\footnotesize 
\[ \mbox{where }\ \opn{C}(f)^m = X^m \oplus {(\bigoplus_{i=m+1}^{m+N-1}}Y^i),\ \ \  
d_{\opn{C}(f)}^m = \left( \begin{array}{c|ccccc}
d&f&0&0&\cdots&0\\
\hline
0&0&1&0&\cdots&0\\
\vdots&\vdots&\vdots&\ddots&\ddots&\vdots\\
0&0&0&\cdots&1&0\\
0&0&0&\cdots&0&1\\
0&-e^{\{N-1\}}&-e^{\{N-2\}}&\cdots& {-e^{\{2\}}}&-e\\
\end{array} \right)
\] }
{\footnotesize 
\[  g^m = \begin{pmatrix}
1\cr 0\cr \vdots\cr 0\cr
\end{pmatrix}, 
\quad
h^m = \begin{pmatrix}
0&1&0&\cdots&0\cr
0&0&1&\ddots&\vdots\cr
\vdots&\vdots&\vdots&\ddots&0\cr
0&0&0&\cdots&1\cr
\end{pmatrix}\ \mbox{ and }\ 
\psi_f ^m = \begin{pmatrix}
f&0&0&\cdots&0\cr
-e&1&0&\cdots&0\cr
0&-e&1&\ddots&\vdots\cr
\vdots&\ddots&\ddots&\ddots&0\cr
0&\cdots&0&-e&1\cr
\end{pmatrix}. 
  \] } 
Thus we have a triangle $Y\xrightarrow{f}X\xrightarrow{g}\opn{C}(f)\xrightarrow{h}\Sigma Y$ in $\KKK_{N} (\BB)$.

\section{Derived category of $N$-complexes}\label{DNcpx}

In this section, we introduce the derived category of $N$-complexes as the Verdier quotient of
the homotopy category with respect to the $N$-quasi-isomorphisms
as in the case of 2-complexes.

\subsection{Homologies of $N$-complexes}\label{HNcpx}

Let $\AA$ be an abelian category, and $\Pj \AA$ (resp., $\Inj \AA$)
the subcategory of $\AA$
consisting of projective (resp., injective) objects of $\AA$.
Let $X$ be an $N$-complex in $\AA$
\[\cdots \to X^{i-1} \xarr{d_X^{i-1}} X^{i} \xarr{d_X^{i}} X^{i+1} \to \cdots.\]
For $0\leq r \leq N$ and $i\in\z$, we define 
\[\begin{aligned}
\opn{Z}^i_{(r)}(X) &:=\opn{Ker}(d_X^{i+r-1}\cdots d_X^i), &
\opn{B}^i_{(r)}(X) &:=\opn{Im}(d_X^{i-1}\cdots d_X^{i-r}), \\
\opn{C}^i_{(r)}(X) &:=\opn{Cok}(d_X^{i-1}\cdots d_X^{i-r}), &
\opn{H}^i_{(r)}(X)&:=\opn{Z}^i_{(r)}(X)/\opn{B}^i_{(N-r)}(X) .
\end{aligned}\]
For example, $\opn{Z}_{(N)}^{n}(X)=\opn{B}^n_{(0)}(X)=X^n$ and $\opn{Z}_{(0)}^{n}(X)=\opn{B}^n_{(N)}(X)=0$ hold. 
With this in mind, using the notation $d^n_{(r)}:= d^n_X|_{\opn{Z}_{(r)}^{n}(X)}$, we can understand a homology as follows
\begin{equation}\label{from Zdiagram}
\opn{H}^n_{(r)}(X)=\opn{Cok}\left(
\opn{Z}^{n-N+r}_{(N)}(X)\xrightarrow{d^{n-N+r}_{(N)}}\cdots
\xrightarrow{d^{n-2}_{(r+2)}}\opn{Z}^{n-1}_{(r+1)}(X)\xrightarrow{d^{n-1}_{(r+1)}}\opn{Z}^n_{(r)}(X)\right).
\end{equation}
For $1\leq r \leq N-1$, we have a pull-back diagram with the canonical inclusion $\iota^n_{(r)}$.
\begin{equation}\label{Ddiagram}
\xymatrix@R1.5em{
0 \ar[r] & \opn{Z}_{(1)}^{n}(X) \ar@{=}[d] \ar[r] & \ar @{} [dr] |{(D_{(r)}^{n})} \opn{Z}_{(r)}^{n}(X) \ar[r]^{d^n_{(r)}} \ar@{^{(}->}[d] ^{\iota^n_{(r)}}
& \opn{Z}_{(r-1)}^{n+1}(X) \ar@{^{(}->}[d]^{\iota^{n+1}_{(r-1)}} \\
0 \ar[r] & \opn{Z}_{(1)}^{n}(X) \ar[r] & \opn{Z}_{(r+1)}^{n}(X) \ar[r]^{d^n_{(r+1)}} & \opn{Z}_{(r)}^{n+1}(X),
}
\end{equation}
Then $(D^{n}_{(r)})$ forms a commutative diagram in Figure 1.

\begin{figure}\label{Zdiagram}
{\tiny
\[
\begin{xy}
(-1,0)*+{X^{n-N+1}}="a1", (21,0)*+{X^{n-N+2}}="a2", 
(41,0)*+{\cdots}="a3", 
(80,0)*+{X^{n-1}}="a5", (100,0)*+{X^{n}}="a6", (120,0)*+{X^{n+1}}="a7", 
(9,-10)*+{Z_{(N-1)}^{n-N+2}}="b1", (20,-10)*+{D_{(N-1)}^{n-N+2}}="b12",
(31,-10)*+{Z_{(N-1)}^{n-N+3}}="b2", 
(51,-10)*+{\cdots}="b3", (70,-10)*+{Z_{(N-1)}^{n-1}}="b4", (80,-10)*+{D_{(N-1)}^{n-1}}="b45", 
(90,-10)*+{Z_{(N-1)}^{n}}="b5", (100,-10)*+{D_{(N-1)}^{n}}="b56", (110,-10)*+{Z_{(N-1)}^{n+1}}="b6", 
(19,-20)*+{Z_{(N-2)}^{n-N+3}}="c1", (30,-20)*+{D_{(N-2)}^{n-N+3}}="c12", 
(41,-20)*+{Z_{(N-2)}^{n-N+4}}="c2",
(60,-20)*+{\cdots}="c3",
(80,-20)*+{Z_{(N-2)}^{n}}="c4", (90,-20)*+{D_{(N-2)}^{n}}="c45", (100,-20)*+{Z_{(N-2)}^{n+1}}="c5", 
(29,-30)*+{Z_{(N-3)}^{n-N+4}}="d1", (51,-30)*+{\ddots}="d2",
(70,-30)*+{\text{\reflectbox{$\ddots$}}}="d3",
(90,-30)*+{Z_{(N-3)}^{n+1}}="d4",
(40,-40)*+{\ddots}="e1", (60,-40)*+{Z_{(2)}^{n}}="e2", (80,-40)*+{\text{\reflectbox{$\ddots$}}}="e3",
(50,-50)*+{Z_{(1)}^{n}}="f1", (60,-50)*+{D_{(1)}^{n}}="f12", (70,-50)*+{Z_{(1)}^{n+1}}="f2", 
(60,-60)*+{0}="g1"
\ar@{->}  "a1";"b1"
\ar@{->}  "b1";"a2"
\ar@{->}  "a2";"b2" 
\ar@{->}  "b2";"a3"
\ar@{->}  "b4";"a5" 
\ar@{->}  "a5";"b5"
\ar@{->}  "b5";"a6" 
\ar@{->}  "a6";"b6"
\ar@{->}  "b6";"a7"
\ar@{->}  "b1";"c1"
\ar@{->}  "c1";"b2" 
\ar@{->}  "b2";"c2"
\ar@{->}  "c2";"b3" 
\ar@{->}  "c3";"b4" 
\ar@{->}  "b4";"c4"
\ar@{->}  "c4";"b5" 
\ar@{->}  "b5";"c5"
\ar@{->}  "c5";"b6"
\ar@{->}  "c1";"d1"
\ar@{->}  "d1";"c2"
\ar@{->}  "c2";"d2"
\ar@{->}  "d3";"c4"
\ar@{->}  "c4";"d4"
\ar@{->}  "d4";"c5"
\ar@{->}  "d1";"e1"
\ar@{->}  "d2";"e2"
\ar@{->}  "e2";"d3"
\ar@{->}  "e3";"d4"
\ar@{->}  "e1";"f1"
\ar@{->}  "f1";"e2"
\ar@{->}  "e2";"f2"
\ar@{->}  "f2";"e3"
\ar@{->}  "f1";"g1"
\ar@{->}  "g1";"f2"
\end{xy}
\]
}
\caption{}
\end{figure}

\begin{definition}\label{acyclic}
We call $X\in\CCC_N(\AA)$ \emph{$N$-acyclic} if $\opn{H}_{(r)}^i(X)=0$ for any $0 < r < N$ and $i \in \mathbb{Z}$. 
\end{definition}

For example, the complex $\mu^i_N(M)$ is $N$-acyclic for any $M\in\AA$ and $i\in\z$. 
An $N$-complex $X$ is $N$-acyclic if and only if there exists some $r$ with $0<r<N$ such that 
$\opn{H}^i_{(r)}(X)=0$ for each integer $i$ \cite{Ka}. 

For $\natural=$nothing$,-, +,\bo$, let $\CCC_{N}^{\natural,\mrm{a}}(\AA)$ (resp., $\KKK_{N}^{\natural,\mrm{a}}(\AA)$) denote the full subcategory of $\CCC_{N}^\natural(\AA)$ (resp., $\KKK^\natural_{N}(\AA)$)
consisting of $N$-acyclic $N$-complexes. 

\begin{proposition}\label{acyclic is triangulated}
We have the following.
\begin{enumerate}
\item $\KKK_{N}^{\natural,\mrm{a}}(\AA)$ is a thick subcategory of 
$\KKK_{N}^{\natural}(\AA)$ for $\natural= -, +, \bo$.
\item $\opn{H}^i _{(r)} (\Sigma X) = \opn{H}^{i +r}_{(N-r)} (X)$ and $\opn{H}^i _{(r)} (\Sigma^{-1} X) = \opn{H}^{i-N+r}_{(N-r)} (X)$ hold for any $X\in\CCC_N(\AA)$.
\end{enumerate}
\end{proposition}

To prove this, 
we recall that $\CCC_N(\AA)$ forms an abelian category.
A sequence $0 \to X \xarr{\alpha} Y \xarr{\beta} Z \to 0$ is exact if and only if 
$0 \to X^i  \xarr{\alpha} Y^i  \xarr{\beta} Z^i  \to 0$ is (not necessarily split) exact in $\AA$ for each $i$.
In this case, for any $0\le r\le N$ and $i\in\z$, we have the following exact sequence \cite{D}.
\begin{equation}\label{ses}
\begin{array}{lllllllll}
\cdots& \xarr{\partial_*}&\opn{H}^{i}_{(r)}(X)\xarr{\alpha_*} \opn{H}^{i}_{(r)}(Y)\xarr{\beta_*} \opn{H}^{i}_{(r)}(Z)
\xarr{\partial_*}\opn{H}^{i+r}_{(N-r)}(X)\xarr{\alpha_*} \opn{H}^{i+r}_{(N-r)}(Y)\xarr{\beta_*} \opn{H}^{i+r}_{(N-r)}(Z)\\
&\xarr{\partial_*}&\opn{H}^{i+N}_{(r)}(X)\xarr{\alpha_*}\opn{H}^{i+N}_{(r)}(Y)\xarr{\beta_*} \opn{H}^{i+N}_{(r)}(Z)
\xarr{\partial_*}\opn{H}^{i+r+N}_{(N-r)}(X)\xarr{\alpha_*} \cdots .&
\end{array}
\end{equation}

\begin{proof}[of Proposition \ref{acyclic is triangulated}]
(2) It is immediate by applying \eqref{ses} to the exact sequences \eqref{hull and cover}.
\par\noindent
(1) It follows from (2) that $\KKK_{N}^{\natural,\mrm{a}}(\AA)$ is closed under $\Sigma$ and $\Sigma^{-1}$.
Let $X\to Y\to Z\to\Sigma X$ be a triangle in $\KKK_N(\AA)$. This comes from a term-wise split short exact sequence. 
Therefore if $X$ and $Y$ belong to $\KKK_N^{\natural,\mrm{a}}(\AA)$, then so does $Z$ by \eqref{ses}.
\end{proof}

As in the classical case, we have the following observation.

\begin{lemma} \label{lem:localqis1}
If $X\in \KKK_{N}^{\mrm{a}}(\AA)$ and $P \in \KKK_{N}^{-}(\cat{Prj}\AA)$
(resp., $I \in \KKK_{N}^{+}(\cat{Inj}\AA)$), then we have
$\opn{Hom}_{\KKK_{N}(\AA)}(P, X) = 0$ (resp., $\opn{Hom}_{\KKK_{N}(\AA)}(X, I) = 0)$.
\end{lemma}

\begin{proof}
Let $f:P\to X$ be as follows.
\[\xymatrix@R1.5em{
\hspace{20pt} P: \ar[d]^{f}\cdots \ar[r] & P^{n-2} \ar[r]^{d_P^{n-2}} \ar[d]^{f^{n-2}}& P^{n-1} \ar[r]^{d_P^{n-1}}  \ar[d]^{f^{n-1}}
& P^{n} \ar[r]  \ar[d]^{f^{n}}& 0 \ar[r]\ar[d]&\cdots\\
\hspace{20pt} X: \cdots \ar[r] & X^{n-2} \ar[r]^{d_X^{n-2}} & X^{n-1} \ar[r]^{d_X^{n-1}} & X^{n} \ar[r]^{d_X^{n}} & X^{n+1} \ar[r] & \cdots. \\
}\]
Since $d_X^{n}f^{n}=0$ and $\opn{H}^n_{(1)}(X)=0$, there is $s^{n}:P^{n} \to X^{n-N+1}$ such that
$f^{n}=d_X^{n-1}\cdots d_X^{n-N+1}s^{n}$.
Since $d_X^{n-1}(f^{n-1} - d_X^{n-2}\cdots d_X^{n-N+1}s^{n}d_P^{n-1})=d_X^{n-1}f^{n-1} - f^{n}d_P^{n-1}=0$,
there is $s^{n-1}:P^{n-1} \to X^{n-N}$ such that
$f^{n-1}=d_X^{n-2}\cdots d_X^{n-N+1}s^{n}d_P^{n-1} + d_X^{n-2}\cdots d_X^{n-N}s^{n-1}$.
Repeating similar argument, we obtain $s^{i}: P^{i} \to X^{i-N+1}$ for $i \leq n$ satisfying \eqref{eq:null homotopic}.
\end{proof}

Now let $\BB$ be an additive category, pick $X\in\CCC_N(\BB)$ and $M\in\BB$. Then we have $N$-complexes  $\Hom_{\BB}(X,M)$ and $\Hom_{\BB}(M,X)$ of abelian groups with $\Hom_{\BB}(M,X)^n:=\Hom_{\BB}(M,X^n)$ and $\Hom_{\BB}(X,M)^n:=\Hom_{\BB}(X^{-n},M)$.
One can easily check the following analogs of \eqref{adjoint} for each $0<r<n$.
{\small\begin{equation}\label{homfrommu1}
\begin{array}{ll}
\Hom_{\CCC_N(\BB)}(\mu_r^{s}(M),X)\simeq\opn{Z}^{s-r+1}_{(r)}(\Hom_{\BB}(M,X)),&\Hom_{\KKK_N(\BB)}(\mu_r^{s}(M),X)\simeq\opn{H}^{s-r+1}_{(r)}(\Hom_{\BB}(M,X)),\\
\Hom_{\CCC_N(\BB)}(X,\mu_r^{s}(M))\simeq\opn{Z}^{-s}_{(r)}(\Hom_{\BB}(X,M)),&
\Hom_{\KKK_N(\BB)}(X,\mu_r^{s}(M))\simeq\opn{H}^{-s}_{(r)}(\Hom_{\BB}(X,M)).
\end{array}
\end{equation}}
We prepare the following observations which will be used later.

\begin{lemma}\label{homfrommu}
Let $X\in\KKK_N(\AA)$, $M\in\AA$, and $0<r<N$ be given. 
 \begin{enumerate}
\item We have a commutative diagram of exact sequences
{\scriptsize
\[\xymatrix@R=1em@C=1em{
&{\Hom}_{\AA}(M,X^{s-N+1})\ar[d]^{d^{\{N-r\}}}\ar[r]&{\Hom}_{\AA}(M,\opn{Z}^{s-r+1}_{(r)}(X))\ar[r]\ar@{=}[d]&{\Hom}_{\KKK_N(\AA)}(\mu_r^{s}(M),X)\ar[r]\ar[d] &0\\
0\ar[r]&{\Hom}_{\AA}(M,\opn{B}^{s-r+1}_{{(N-r)}}(X))\ar[r]&{\Hom}_{\AA}(M,\opn{Z}^{s-r+1}_{(r)}(X))\ar[r]&{\Hom}_{\AA}(M,\opn{H}^{s-r+1}_{(r)}(X))\ar[r]&\Ext^1_{\AA}(M,\opn{B}^{s-r+1}_{{(N-r)}}(X))}
\]
}
\item If $M$ is projective in $\AA$, then $\Hom_{\KKK_{N}(\AA)}(\mu_{r}^{s}(M), X) \simeq \Hom_{\AA}(M,\opn{H}_{(r)}^{s-r+1}(X))$.
\item If $X\in\KKK_N(\Inj\AA)$ is $N$-acyclic, then $\Hom_{\KKK_{N}(\AA)}(\mu^{s}_{r}(M),X)\simeq\Ext_{\AA}^{1}(M,\opn{Z}^{s-N+1}_{(N-r)}(X))$.
\end{enumerate}
\end{lemma}

\begin{proof}
(1) The upper sequence is exact by \eqref{homfrommu1} and $\opn{Z}^{s-r+1}_{(r)}(\Hom_{\AA}(M,X))\simeq\Hom_{\AA}(M,\opn{Z}^{s-r+1}_{(r)}(X))$. 
The lower one is clearly exact.
\par\noindent
(2) Immediate from (1).
\par\noindent
(3) We have a short exact sequence $0\to\opn{Z}^{s-N+1}_{(N-r)}(X)\to X^{s-N+1}\to\opn{Z}^{s-r+1}_{(r)}(X)\to0$.
Applying $\Hom_{\AA}(M,-)$ and using injectivity of $X^{s-N+1}$, we have an exact sequence
\[\Hom_{\AA}(M,X^{s-N+1})\to\Hom_{\AA}(M,\opn{Z}^{s-r+1}_{(r)}(X))\to\Ext^1_{\AA}(M,\opn{Z}^{s-N+1}_{(N-r)}(X))\to0.\]
Comparing with the upper exact sequence in (1), we have the desired isomorphism.
\end{proof}

\begin{lemma}\label{0cpx}
For a commutative diagram \eqref{Ddiagram}, the following hold. 
\begin{enumerate}
\item If $\opn{H}^{n} _{(r)} (X) =0$, then $(D^{n+r-1-s}_{(s)})$ is an exact square for any $r\le s\le N-1$. 
In particular, $(D^{n-1}_{(r)} + D^{n-2}_{(r+1)} + \cdots + D^{n-N+r}_{(N-1)})$ is an exact square. 
\item $X$ is $N$-acyclic if and only if $d^{n}_{(r+1)}$ is an epimorphism for any $0 < r < N$ and $n\in \mathbb{Z}$.
\item $X$ is isomorphic to $0$ in $\KK{A}$ if and only if $d^{n}_{(r+1)}$ is a split epimorphism 
for any $0 < r < N$ and $n \in \mathbb{Z}$.
\end{enumerate}
\end{lemma}

\begin{proof}
(1) (2) The assertions immediately follow  from \eqref{from Zdiagram}.
\par\noindent
(3) We prove the `only if' part. 
Clearly $d^{n}_{(r+1)}$ is a split epimorphism for $X=\mu^s_N(M)$.
Since every projective-injective object of $\CCC_N(\AA)$ is in  
$\Add\{\mu^s_N(M)\mid s\in\z,\ M\in\AA\}$,
the assertion follows.

To show the converse, set 
$\opn{W}^n_{(r)}(X):=\bigoplus_{i=0}^{r-1}\opn{Z}^{n+i}_{(1)}(X)$ for $1\le r \le N$. 
Then we have natural morphisms $p^n_{(r)}:=\left(\begin{smallmatrix}0&1\end{smallmatrix}\right):\opn{W}^{n}_{(r)}(X)\to \opn{W}^{n+1}_{(r-1)}(X)$
and $i^n_{(r)}:=\left(\begin{smallmatrix}1\\ 0\end{smallmatrix}\right):\opn{W}^{n}_{(r)}(X)\to \opn{W}^{n}_{(r+1)}(X)$.
We show the existence of an isomorphism $a^n_{(r)}:\opn{W}^{n}_{(r)}(X)\to \opn{Z}^n_{(r)}(X)$ such that the following diagram commute.
{\small\[\xymatrix@R=.3em@C=5em{
\opn{W}^n_{(r)}(X)\ar[rr]|{p^n_{(r)}}\ar[dr]|{a^n_{(r)}}\ar[dd]|{i^n_{(r)}}&&\opn{W}^{n+1}_{(r-1)}(X)\ar[dd]|(.4){i^{n+1}_{(r-1)}}\ar[dr]|{a^{n+1}_{(r-1)}}\\
&\opn{Z}^n_{(r)}(X)\ar[rr]|(.3){d^n_{(r)}}\ar[dd]|(.4){\iota^n_{(r)}}&&\opn{Z}^{n+1}_{(r-1)}(X)\ar[dd]|{\iota^{n+1}_{(r-1)}}\\
\opn{W}^n_{(r+1)}(X)\ar[rr]|(.3){p^n_{(r+1)}}\ar[dr]|{a^n_{(r+1)}}&&\opn{W}^{n+1}_{(r)}(X)\ar[dr]|{a^{n+1}_{(r)}}\\
&\opn{Z}^n_{(r+1)}(X)\ar[rr]|{d^n_{(r+1)}}&&\opn{Z}^{n+1}_{(r)}(X).
}\]}
For $r=1$, set $a^n_{(1)}=1$.
Suppose $r>1$ and that we have defined $a^n_{(i)}$ for any $n,i \in \mathbb{Z}$ with $0<i \le r$. 
Applying Lemma \ref{split square} to the exact square 
{\small \[ \xymatrix@R=1.5em@C=5em{
\opn{W}^n_{(r)}(X)=\opn{Z}^n_{(1)}(X) \oplus \opn{W}^{n+1}_{(r-1)}(X)\ar[rr]^{p^n_{(r)}} \ar[d]^{\iota ^n_{(r)}  a^n_{(r)} } &&\opn{W}^{n+1}_{(r-1)}(X) \ar[d]^{i^{n+1}_{(r-1)}} \\
\opn{Z}^n_{(r+1)}(X) \ar[rr]^{ ( a^{n+1}_{(r)} )^{-1}  d^n_{(r+1)}} &&\opn{W}^{n+1}_{(r-1)}(X)\oplus \opn{Z}^{n+r}_{(1)}(X) =\opn{W}^{n+1}_{(r)}(X),} \]}
we get an isomorphism 
$a^n_{(r+1)}:\opn{W}^{n}_{(r+1)}(X)\to \opn{Z}^n_{(r+1)}(X)
$ as desired.
\par\noindent
Consequently we have an isomorphism $a^n_{(N)}:\opn{W}^n_{(N)}(X)=\bigoplus_{i=0}^{N-1}\opn{Z}^{n+i}_{(1)}(X)\to \opn{Z}^n_{(N)}(X)=X^n$.
Since $d^n=\iota^{n+1}_{(N-1)}d^n_{(N)}:X^n\to X^{n+1}$ holds, it is easy to check $X\simeq\bigoplus_{n\in\z}\mu^n_N(\opn{Z}^n_{(1)}(X))$ in $\CCC_N(\AA)$. Thus $X$ is zero in $\KK{A}$.
\end{proof}

\begin{definition}\label{N-qis}
A morphism $f: X \to Y$ of $\KKK_{N}(\AA)$ is called an \emph{$N$-quasi-isomorphism}
if $\opn{H}^i_{(r)}(f): \opn{H}^i_{(r)}(X) \to\opn{H}^i_{(r)}(Y)$
is an isomorphism for any $0 < r < N$ and $i \in \mathbb{Z}$, or equivalently by \eqref{ses}, the mapping cone $\opn{C}(f)$ is $N$-acyclic.
For $\natural=$nothing$,+,-,\bo$,
the \emph{derived category} of $N$-complexes is defined as the quotient category
\[\DDD_{N}^{\natural}(\AA)=\KKK_{N}^{\natural}(\AA)/\KKK_{N}^{\natural,\mrm{a}}(\AA).\]
\end{definition}

By definition, a morphism in $\KKK_{N}^{\natural}(\AA)$ is an $N$-quasi-isomorphism 
if and only if it is an isomorphism in $\DDD_{N}^{\natural}(\AA)$.

\begin{proposition} \label{prop:dertrig}
\begin{enumerate}
\item If $0 \arr X \xarr{f} Y \xarr{g} Z \arr 0$ is an exact sequence in the abelian category $\CCC_{N}(\AA)$, then it can be embedded into a triangle $X \xarr{f} Y \xarr{g} Z \xarr{h} \Sigma X$ in $\DDD_{N}(\AA)$.
\item For any triangle $X \xarr{f} Y \xarr{g} Z \xarr{h} \Sigma Y$ in $\DD{A}$, we have a long exact sequence
\small{
\[
\begin{array}{lllllllll}
\cdots& \xarr{h_*}&\opn{H}^{i}_{(r)}(X)\xarr{f_*} \opn{H}^{i}_{(r)}(Y)\xarr{g_*} \opn{H}^{i}_{(r)}(Z)
\xarr{h_*}\opn{H}^{i+r}_{(N-r)}(X)\xarr{f_*} \opn{H}^{i+r}_{(N-r)}(Y)\xarr{g_*} \opn{H}^{i+r}_{(N-r)}(Z)\\
&\xarr{h_*}&\opn{H}^{i+N}_{(r)}(X)\xarr{f_*}\opn{H}^{i+N}_{(r)}(Y)\xarr{g_*} \opn{H}^{i+N}_{(r)}(Z)
\xarr{h_*}\opn{H}^{i+r+N}_{(N-r)}(X)\xarr{f_*} \cdots .&
\end{array}
\]
}
\end{enumerate}
\end{proposition}

\begin{proof}
(1) We have the following commutative diagram of exact sequences in $\CCC_{N}(\AA)$.
{\small\[\xymatrix@R=1em{
& 0 \ar[d] & 0 \ar[d] \\
0 \ar[r] & X \ar[r] \ar[d]^{f} & I(X) \ar[r] \ar[d]^{\psi_f} & \Sigma X \ar[r] \ar@{=}[d] & 0 \\
0 \ar[r] & Y \ar[r]^{u}\ar[d]^{g} & \opn{C}(f) \ar[r]^{v}\ar[d]^{s} & \Sigma X \ar[r] & 0 \\
& Z \ar@{=}[r] \ar[d] & Z \ar[d] \\
&0  & 0
}\]}
Then $X \xarr{f} Y \xarr{u} \opn{C}(f)\xarr{v} \Sigma X$ is a triangle in $\KKK_{N}(\AA)$.
Since $I(X)$ is $N$-acyclic, $s$ is an $N$-quasi-isomorphism.
Thus we have a triangle $X \xarr{f} Y \xarr{su=g} Z\xarr{vs^{-1}} \Sigma X$ in $\DDD_{N}(\AA)$.
\par\noindent
(2) We have only to verify the assertion for the triangle $X \xarr{f} Y \xarr{} \opn{C}(f) \xarr{} \Sigma Y$.
Applying \eqref{ses} to a short exact sequence $0\to X\to Y\oplus I(X)\to \opn{C}(f)\to0$ in $\CCC_N(\AA)$, we get the desired sequence.
\end{proof}

\begin{definition}[Truncations] \label{dfn:trunc}
For an $N$-complex $X = (X^{i}, d^{i})$, set
\[ 
{\sigma}_{\leq n}X  : \cdots \xarr{d^{n-N}} X^{n-N+1} \xarr{d_{(N)}^{n-N+1}}  \opn{Z}_{(N-1)}^{n-N+2}(X) 
\xarr{d_{(N-1)}^{n-N+2}} \cdots \xarr{d_{(2)}^{n+1}} \opn{Z}_{(1)}^{n}(X)  \arr 0 \arr \cdots.\\ 
\]
\end{definition}

\begin{lemma}\label{trunc01}
For an $N$-complex $X = (X^{i}, d^{i})$ and an integer $n$, the following hold. 
\begin{enumerate}
\item
$\opn{H}_{(r)}^{i}(\sigma_{\leq n}(X)) \simeq \opn{H}_{(r)}^{i}(X)$ for any $0 < r < N$ and $i +r \leq n+1$.
\item If $\opn{H}_{(r)}^i(X)=0$ holds for any $0 < r < N$ and $i\geq n+1$,
then the canonical injection $\sigma_{\leq n}X \to X$ is an $N$-quasi-isomorphism.
\end{enumerate}
\end{lemma}

\begin{proof}
(1) If $i+r\leq n+1$, then $\opn{Z}_{(r)}^{i}(X)$ is the kernel of 
${d^{\{r\}}}:  \opn{Z}_{(n-i+1)}^{i} (X) \to X^{i+r}$ 
which maps into $\opn{Z}_{(n-i-r+1)}^{i+r}(X)$. Hence $\opn{Z}_{(r)}^{i}(X)=\opn{Z}_{(r)}^{i}(\sigma_{\leq n}X) $. 
Clearly $\opn{B}_{(N-r)}^{i}(\sigma_{\leq n}X) =\opn{B}_{(N-r)}^{i}(X)$.
\par\noindent
(2) It remains to show $\opn{H}_{(r)}^{i}(\sigma_{\leq n}(X)) \simeq \opn{H}_{(r)}^{i}(X)$ for $i \le n$ and $i+r>n+1$.
Since $\opn{Z}_{(r)}^{i}(\sigma_{\leq n}X) =\opn{Z}_{(n-i+1)}^{i}(X)$ holds,
we have a commutative diagram
\[\xymatrix@R=1em{
\opn{Z}_{(n-i+N-r+1)}^{i-N+r}(X) \ar[r]^{\quad d^{\{N-r\}}} \ar@{^{(}->}[d]& \opn{Z}_{(n-i+1)}^{i}(X)\ar@{^{(}->}[d] \ar[r]
&  \opn{H}^{i}_{(r)}(\sigma_{\leq n}X) \ar[r] & 0 \\
X^{i-N+r} \ar[r]^{d^{\{N-r\}}}  & \opn{Z}_{(r)}^{i}(X) \ar[r] & \opn{H}^{i}_{(r)}(X) \ar[r] & 0
}\]
of exact sequences. 
The left square is exact. Indeed it follows from Lemmas \ref{prop:pullpush3} and \ref{pullback2} 
since $(D^j_{(s)} )$ is an exact square for $j+s \geq n+1$ by Lemma \ref{0cpx}(1). 
Thus we have the desired isomorphism.
\end{proof}

\begin{proposition}\label{b+-}
Let $\natural = +, -, \bo$.
The canonical functors $\DDD_{N}^{\bo} (\AA)\to\DDD_{N}^{\natural} (\AA)\to\DDD_{N}(\AA) $ are fully faithful.
Therefore $\DDD_{N}^{\natural}(\AA)$ is equivalent to the full subcategory of $\DDD_{N}(\AA)$ 
consisting of objects in $\KKK_{N}^{\natural}(\AA)$.
\end{proposition}

\begin{proof} 
We only show that $\DDD_{N}^{-} (\AA)\to \DDD_{N}(\AA)$ is fully faithful.
Let $f:X\to Y$ be any morphism with $X\in\KKK_{N}^{-} (\AA)$ and $Y\in \KKK_{N}^{\mrm{a}}(\AA)$.
For sufficiently large $n$, $f$ factors through the natural morphism $\sigma_{\le n}(Y)\to Y$.
Since $\sigma_{\le n}(Y)$ belongs to $\KKK_{N}^{-, \mrm{a}} (\AA)$ by Lemma \ref{trunc01}(2), 
we get the conclusion from Lemma \ref{fully faithful}.
\end{proof}

\subsection{Elementary morphisms}\label{EMN}

In this subsection, we introduce the $N$-complex version of {\it an elementary map of degree $i$} in the sense of Verdier \cite{Ve}.
We start with the following observation.

\begin{definition-proposition}\label{elmap02}
For an object $X : \cdots \to X^{i-1} \xarr{d_X^{i-1}} X^{i} \xarr{d_X^{i}} X^{i+1} \to \cdots$ in $\CCC_{N}(\AA)$ and 
a morphism $u: M \to X^i$ in $\AA$,  we take successive pull-backs
\[\xymatrix@R1em{
\ar @{} [drr] |{(E^{i-r-1})} Y^{i-r-1} \ar[d]_{u^{i-r-1}} \ar[rr]^{d'^{i-r-1}} && Y^{i-r} \ar[d]^{u^{i-r}} \\
X^{i-r-1} \ar[rr]_{d_X^{i-r-1}} && X^{i-r}
}\]
for $0 \leq r < N-1$, where $Y^i=M$ and $u^i=u$. Then there are a morphism $d'^{i-N}:X^{i-N} \to Y^{i-N+1}$ in $\AA$ and a morphism
{\small
\[\xymatrix@R1.5em{
\opn{V}_{i}(X,u): \ar[d]^{p_i(u)}&
\cdots \ar[r] &X^{i-N} \ar@{=}[d] \ar[r]^{d'^{i-N}} &  \ar @{} [dr] |{(E^{i-N+1})}Y^{i-N+1} \ar[d]_{u^{i-N+1}} \ar[r]^{\quad d'^{i-N+1}} 
& \cdots \ar[r] & \ar @{} [dr] |{(E^{i-1})} Y^{i-1}  \ar[d]_{u^{i-1}} \ar[r]^{d'^{i-1}} & M  \ar[d]^{u} \ar[r]^{d_X^{i}u}  & X^{i+1}  \ar@{=}[d] \ar[r] & \cdots \\
X: & \cdots \ar[r] &X^{i-N} \ar[r]_{d_X^{i-N}} & X^{i-N+1} \ar[r]_{\quad d_X^{i-N+1}} 
& \cdots \ar[r] & X^{i-1} \ar[r]_{d_X^{i-1}} & X^i \ar[r]_{d_X^{i}} & X^{i+1} \ar[r] & \cdots
}\]
}
in $\CCC_{N}(\AA)$. Moreover the following conditions are equivalent.
\begin{enumerate}
\item  $p_i(u)$ is an $N$-quasi-isomorphism.
\item  The commutative diagram $(E^{i-N+1}+ \cdots +E^{i-1})$ is an exact square.
\item  The commutative diagrams $(E^{i-N+1}), \cdots , (E^{i-1})$ are exact squares.
\item  $(u\ d^{\{N-1\}}):M\oplus X^{i-N+1}\to X^i$ is an epimorphism.
\end{enumerate}
\end{definition-proposition}

\begin{proof}
Set $Y=V_i(X,u)$ and $\tilde{u}=p_i(u)$.
\par\noindent
(2) $\Leftrightarrow$ (3)$\Leftrightarrow$ (4). These are clear from Lemmas~\ref{pullback2} and \ref{prop:pullpush3}.
\par\noindent
(1) $\Rightarrow$ (4). 
The morphism $\tilde{u}$ induces a morphism $\overline{u}: \overline{Y} \to\overline{X}$ of $2$-complexes as follows: 
{\small 
\[ \xymatrix@R1.5em{
\overline{Y}: \ar[d]^{\overline{u}} & \ar[r]^{d_Y ^{\{ N-1 \} } } & 
Y^{i-N} \ar[r]^{d_Y} \ar@{=}[d]  &\ar @{} [dr] |{(E)}  Y^{i-N+1} \ar[d] _{u^{i-N+1} } \ar[r] ^(.6){~d_Y ^{\{ N-1 \} } } & 
M \ar[r]^{d_Y} \ar[d] ^{u} & Y^{i+1} \ar@{=}[d] \ar[r] ^{d_Y ^{\{ N-1 \} } } &
Y^{i+N} \ar[r]^{d_Y} \ar@{=}[d] & \\
\overline{X}: &\ar[r]_{d_X ^{\{ N-1 \} } } & 
X^{i-N} \ar[r]_{d_X}  & X^{i-N+1} \ar[r] _(.6){~d_X ^{\{ N-1 \} } } & 
X^{i} \ar[r]_{d_X} & X^{i+1}  \ar[r] _{d_X ^{\{ N-1 \} } } &
X^{i+N} \ar[r]_{d_X}  &\\
} \] }
The assumption forces $\overline{u}$ to be a $2$-quasi-isomorphism.  
Then \cite[III. 2.1.2(c]{Ve} implies that 
$(u\ d^{\{N-1\}}):M\oplus X^{i-N+1}\to X^i$ is an epimorphism.
\par\noindent
(3) $\Rightarrow$ (1). 
We shall show that 
$\opn{H}_{(r)}^{i-s}(\tilde{u}) : \opn{H}_{(r)}^{i-s}(Y) \to \opn{H}_{(r)}^{i-s}(X)$ is an isomorphism for 
each $s \in \mathbb{Z}$, $0< r<N$. 
Set the commutative squares (A), (B), (C), (D) as follows: 
\[ \xymatrix@R1.5em{ 
\ar @{} [dr] |{(A)}  Y^{i-N-s} \ar[d]_{} \ar[r]^{d_Y^{\{ r\}}} &  \ar @{} [dr] |{(B)} Y^{i-N-s+r}  \ar[r]^(0.6){d_Y^{\{N- r\}} } \ar[d]^{} & 
\ar @{} [dr] |{(C)} Y^{i-s}  \ar[r]^{d_Y^{\{ r\}}} \ar[d]&  \ar @{} [dr] |{(D)} Y^{i-s+r}  \ar[r]^{d_Y^{\{N- r\}}} \ar[d]& 
Y^{i+N-s}  \ar[d]\\
X^{i-N-s}  \ar[r]_{d_X^{\{ r\}}} &  X^{i-N-s+r}  \ar[r]_{d_X^{\{N- r\}} }& 
X^{i-s}  \ar[r]_{d_X^{\{ r\}}} &  X^{i-s+r}  \ar[r]_{d_X^{\{N- r\}}} & 
X^{i+N-s}  \\} \]

Assume that (A) and (C) are exact. Consider the diagram with exact rows 
\[ \xymatrix@R1.5em{ 
\opn{C}^{i-N-s+r} _{(r)} (Y) \ar[d]^{\opn{C}^{i-N-s+r}_{(r)}(\tilde{u})} \ar[r] & 
\opn{Z}^{i-s}_{(r)} (Y) \ar[d]^{\opn{Z}^{i-s}_{(r)}(\tilde{u})} \ar[r] & 
\opn{H}^{i-s}_{(r)} (Y) \ar[d]^{\opn{H}^{i-s}_{(r)}(\tilde{u})} \ar[r] &
0\\
\opn{C}^{i-N-s+r} _{(r)} (X)  \ar[r]& 
\opn{Z}^{i-s}_{(r)} (X) \ar[r] & 
\opn{H}^{i-s}_{(r)} (X)  \ar[r] &
0.\\} \]
Lemma~\ref{prop:pullpush3} implies that $\opn{C}^{i-N-s+r}_{(r)}(\tilde{u})$ and $\opn{Z}^{i-s}_{(r)}(\tilde{u})$ are isomorphisms. 
Hence so is $\opn{H}^{i-s}_{(r)}(\tilde{u})$. 
Similarly $\opn{H}^{i-s}_{(r)}(\tilde{u})$ is an isomorphism provided that (B) and (D) are exact.
Therefore it is enough to show that either (A), (C) or (B), (D) are exact. 
To prove this, notice that for any integer $j$ other than $i-N$ or $i$, the following square is exact. 
\[ \xymatrix@R1em{ Y^j \ar[r]^{d_Y^j} \ar[d]_{u^j} & 
Y^{j+1} \ar[d]^{u^{j+1}} \\
X^j \ar[r]^{d_X^j} & 
X^{j+1} \\ } \]
Lemma \ref{pullback2} (1)$\Rightarrow$(2) implies that (B) and (D) are exact  if $s\in\{0,1,\ldots,r-1\}$, 
otherwise (A) and (C) are exact.
Therefore one of the above two conditions holds.
\end{proof}

\subsection{Resolutions of $N$-complexes}

The aim of this subsection is to establish Theorems \ref{cor:subeqv}, \ref{cor:subeqv03}
which are well-known for the classical case $N=2$.

For a full additive subcategory $\BB$ of an abelian category $\AA$ and $\natural =$nothing$,-, +, \bo$,
we denote by $\CCC_{N}^{\natural,\mrm{a}}(\BB)$ (resp., $\CCC_{N}^{\natural, \bo}(\BB)$, 
$\CCC_{N}^{\natural,-}(\BB)$, $\CCC_{N}^{\natural,+}(\BB)$)
the full subcategory of $\CCC_{N}^{\natural}(\BB)$ 
consisting of $N$-complexes $X$  satisfying that $\opn{H}_{(r)}^{i}(X) = 0$ for any $0 < r < N$ and 
for all (resp. all but finitely many, sufficiently large, sufficiently small) 
$i\in \mathbb{Z}$.
The corresponding subcategory of $\KKK_{N}^{{\natural}}(\BB )$ is denoted by 
 $\KKK_{N}^{\natural,\mrm{a}}(\BB)$ 
 (resp., $\KKK_{N}^{\natural, \bo}(\BB)$, $\KKK_{N}^{\natural,-}(\BB)$, $\KKK_{N}^{\natural,+}(\BB)$).

\begin{theorem}\label{cor:subeqv}
The following hold for $\natural$=nothing$,\bo$.
\begin{enumerate}
\item If $\AA$ has enough projectives, then $(\KKK_{N}^{-,\natural}(\cat{Prj}\AA),\KKK_{N}^{-,\mrm{a}}(\AA))$ is a stable t-structure
in $\KKK_{N}^{-,\natural}(\AA)$ and we have triangle equivalences
$\KKK_{N}^{-}(\cat{Prj}\AA) \simeq \DDD_{N}^{-}(\AA)$ and $\KKK_{N}^{-,\bo}(\cat{Prj}\AA) \simeq \DDD_{N}^{\bo}(\AA)$.
\item If $\AA$ has enough injectives, then $(\KKK_{N}^{+,\mrm{a}}(\AA),\KKK_{N}^{+,\natural}(\cat{Inj}\AA))$ is a stable t-structure
in $\KKK_{N}^{+,\natural}(\AA)$ and we have triangle equivalences
$\KKK_{N}^{+}(\cat{Inj}\AA) \simeq \DDD_{N}^{+}(\AA)$ and $\KKK_{N}^{+,\bo}(\cat{Inj}\AA) \simeq \DDD_{N}^{\bo}(\AA)$.
\end{enumerate}
\end{theorem}

Our proof of Theorem \ref{cor:subeqv} is based on Verdier's method \cite[III, Section 2.2]{Ve}.

\begin{definition}\label{Ve01}
Let $\MM$ be an additive full subcategory of $\AA$ satisfying the following.
\begin{enumerate}
\item[$V_1$]  For any epimorphism $u:X \to L$ with $X\in \AA$ and $L \in \MM$, there is an epimorphism
$v : L' \to L$ with $L' \in \MM$ which factors through $u$.
\item[$V_2$]  For any exact sequence $0 \to X \to L_n \to \cdots \to L_0 \to 0$ with $L_0, \cdots , L_n \in \MM$,
there is an epimorphism $L' \to X$ with $L' \in \MM$.
\end{enumerate}
Let $\widehat{\MM}$ be the full subcategory of $\AA$ consisting of objects $X$ satisfying the following conditions.
\begin{enumerate}
\item $X$ has an \emph{$\infty$-$\MM$-presentation}, that is, an exact sequence $\cdots \to L_n \to \cdots \to L_1 \to L_0 \to X \to 0$ with $L_i \in \MM$ for any $i \geq 0$.
\item  For any exact sequence $0 \to X' \to L_n \to \cdots \to L_0 \to X \to 0$ with $L_0, \cdots , L_n \in \MM$, 
$X'$ has an $\infty$-$\MM$-presentation.
\end{enumerate}
Obviously we have $\MM \subset \widehat{\MM}$. 
\end{definition}

For example, $\MM=\Proj\AA$ satisfies ($V_1$) and ($V_2$). If $\AA$ has enough projectives, then $\widehat{\MM}=\AA$. 
\begin{lemma}\label{set01}
Let $\MM$ be an additive full subcategory of $\AA$ satisfying ($V_1$) and ($V_2$). 
\begin{enumerate}
\item \cite[III 2.2.4]{Ve} For an exact sequence $0 \to X \to Y \to Z \to 0$, if two out of three terms belong to 
$\widehat{\MM}$, then so does the other.
\item 
For an epimorphism $\rho: X \to L$ with $L \in \widehat{\MM}$, there exists 
a morphism $\mu : M \to X$ with $M \in \MM$ such that 
$\rho \mu$ is an epimorphism. 
\end{enumerate}
 \end{lemma}
 
 \begin{proof} 
 (2) Take an epimorphism $\pi : M_0 \to L$ with $M_0 \in \MM$ and a pull-back diagram 
 \[\xymatrix@R1em{
K\ \ar[r]^{\rho '} \ar[d]_{\pi '} & M_0 \ar[d]^{\pi } \\
X\ \ar[r]^{\rho} & L.}\] 
Then $\rho '$ and $\pi '$ are epimorphisms. 
The condition ($V_1$) gives a morphism 
$k: M \to K$ with $M \in \MM$ such that 
$\rho ' k$ is an epimorphism. 
Set $\mu = \pi ' k$, then $\rho \mu = \pi \rho ' k$ is an epimorphism. 
 \end{proof}

\begin{proposition} \label{prop:localqis2}
Under the conditions ($V_1$) and ($V_2$), we have the following.
\begin{enumerate}
\item Given $X \in \CCC_{N}^{-}(\widehat{\MM})$,  there exists an $N$-quasi-isomorphism $s:L\to X$ 
with $L \in \CCC_{N}^{-}(\MM)$.
\item 
We have $\KKK^{-,\natural}_N(\widehat{\MM})
= \KKK^{-,\natural}_N(\MM)* \KKK^{-,\mrm{a}}_N(\widehat{\MM})$ for $\natural=$nothing$,\bo$.
\item We have a stable t-structure $\left(\frac{\KKK^{-,\natural}_N(\MM)}{\KKK^{-,\mrm{a}}_N(\MM)},\frac{\KKK^{-,\mrm{a}}_N(\widehat{\MM})}{\KKK^{-,\mrm{a}}_N(\MM)}\right)$ in 
$\frac{\KKK^{-,\natural}_N(\widehat{\MM})}{\KKK^{-,\mrm{a}}_N(\MM)}$ and a triangle equivalence
$\frac{\KKK^{-,\natural}_N(\MM)}{\KKK^{-,\mrm{a}}_N(\MM)}\simeq
\frac{\KKK^{-,\natural}_N(\widehat{\MM})}{\KKK^{-,\mrm{a}}_N(\widehat{\MM})}$ for $\natural=$nothing$,\bo$.
\end{enumerate}
\end{proposition}

\begin{proof}
(1) We shall construct a series of $N$-quasi-isomorphisms $v_{n+1}:  L_{n} \to L_{n+1}$ satisfying 
$L_n \in \CCC^{-}_N (\widehat{\MM})$, $L_n^i \in \MM$ ($i>n$)
and $v_{n+1}^i=\mathrm{id}$ ($i>n+1$) by an induction on $n$.
\par\noindent
We set $L_m =X$ and $v_m=\mathrm{id}_X$ for $m$ large enough. 
Suppose we get $L_{n}$ and $v_{n+1}$.
Since $L_{n}^n \in \widehat{\MM}$, there exists an epimorphism $f: M \to L_{n}^n$ with $M \in \MM$. 
Then $L_{n-1} = V_{n} ( L_{n} , f )$ and $v_{n} = p_{n} (f)$ satisfy the conditions above. 
Indeed, $v_{n}$ is an $N$-quasi-isomorphism by Definition-Proposition \ref{elmap02}(4)$\Rightarrow$(1), $L_{n-1}^i \in \MM$ $(i> n-1)$ and $v_{n}^i = \mathrm{id}$ ($i>n$) by the construction, and 
$L_{n-1}^{i} \in \widehat{\MM}$ ($i \leq n-1$) by Lemma \ref{set01}(1).
\par\noindent
Since $v_{n+1}^i: L_{n}^i \to L_{n+1}^i$ ($i>n+1$) is an identity, the canonical morphism $L:=\underset{\leftarrow}{\lim}L_n \to X$ gives a desired $N$-quasi-isomorphism.
\par\noindent
(2) It suffices to prove "$\subset$".
Given an object $X \in \KKK_{N}^{-}(\widehat{\MM})$, there exists an $N$-quasi-isomorphism 
$L \stackrel{s}{\to} X$ with $L \in \KKK_{N}^{-}(\MM)$ by (1).
Then $\mrm{C}(s)\in\KKK_{N}^{-}(\widehat{\MM})$ is $N$-acyclic, and we have
$\KKK^{-}_N(\widehat{\MM})\subset \KKK^{-}_N(\MM)* \KKK^{-,\mrm{a}}_N(\widehat{\MM})$.  
If $X \in \KKK_{N}^{-, \bo}(\widehat{\MM})$, then $L \in \KKK_{N}^{-, \bo}(\MM)$ holds, and hence
 $\KKK^{-, \bo}_N(\widehat{\MM})\subset \KKK^{-, \bo}_N(\MM)*\KKK^{-,\mrm{a}}_N(\widehat{\MM})$. 
\par\noindent
(3) Set $\UU=\KKK^{-,\natural}_N(\MM)$ and $\VV=\KKK^{-,\mrm{a}}_N(\widehat{\MM})$. Then $\UU*\VV=\KKK^{-,\natural}_N(\widehat{\MM})$ holds by (2).
Applying Lemma \ref{fully faithful}, we have a stable t-structure
$(\frac{\UU}{\UU\cap\VV},\frac{\VV}{\UU\cap\VV})=\left(\frac{\KKK^{-,\natural}_N(\MM)}{\KKK^{-,\mrm{a}}_N(\MM)},\frac{\KKK^{-,\mrm{a}}_N(\widehat{\MM})}{\KKK^{-,\mrm{a}}_N(\MM)}\right)$ in 
$\frac{\UU*\VV}{\UU\cap\VV}=\frac{\KKK^{-,\natural}_N(\widehat{\MM})}{\KKK^{-,\mrm{a}}_N(\MM)}$ and triangle equivalences
$\frac{\KKK^{-,\natural}_N(\MM)}{\KKK^{-,\mrm{a}}_N(\MM)}\simeq
\frac{\UU}{\UU\cap\VV}\simeq\frac{\UU*\VV}{\VV}=\frac{\KKK^{-,\natural}_N(\widehat{\MM})}{\KKK^{-,\mrm{a}}_N(\widehat{\MM})}$.
\end{proof}

\begin{proof}[of Theorem \ref{cor:subeqv}]
We only prove (1) since (2) is the dual. 
Set $\MM=\cat{Prj}\AA$, then $\widehat{\MM}=\AA$.
By Lemma \ref{lem:localqis1}, we have $\KKK^{-,\mrm{a}}_N(\MM)=0$.
By Proposition \ref{prop:localqis2}(3), we have a stable t-structure $(\KKK_{N}^{-,\natural}(\cat{Prj}\AA),\KKK_{N}^{-,\mrm{a}}(\AA))$ in $\KKK_{N}^{-,\natural}(\AA)$
and a triangle equivalence $\KKK^{-,\natural}_N(\cat{Prj}\AA)\simeq\frac{\KKK^{-,\natural}_N(\AA)}{\KKK^{-,\mrm{a}}_N(\AA)}$. 
This is $\DDD^-_N(\AA)$ if $\natural$=nothing, and $\DDD^{\bo}_N(\AA)$ if $\natural=\bo$ by Proposition \ref{b+-}.
\end{proof}
 
Recall that an abelian category $\AA$ is an \emph{$Ab3$-category} (resp., \emph{$Ab3^*$-category}) provided that
it has an arbitrary coproduct (resp., product) of objects. 
It is clear that coproducts (resp., products) preserve cokernels (resp., kernels).
Moreover $\AA$ is an \emph{$Ab4$-category} (resp., \emph{$Ab4^*$-category}) provided that
it is an $Ab3$-category (resp., $Ab3^*$-category), and that
the coproduct (resp., product) of monomorphisms (resp., epimorphisms) is monic (resp., epic)
(see e.g. \cite{Po}).

\begin{definition}[cf. \cite{BN,Sp}]\label{dfn:spcpx}
We say that 
$X\in\KKK_N(\AA)$ is \emph{$\KKK$-projective} if $\opn{Hom}_{\KKK_{N}(\AA)}(X, \\ \KKK_{N}^{\mrm{a}}(\AA))=0$.
We say that $X\in\KKK_N(\AA)$ is \emph{$\KKK$-injective} if $\opn{Hom}_{\KKK_{N}(\AA)}(\KKK_{N}^{\mrm{a}}(\AA),X)=0$.
We denote by $\KKK_{N}^{\mrm{p}}(\AA)$ (resp., $\KKK_{N}^{\mrm{i}}(\AA)$)
the full triangulated subcategory of $\KKK_{N}(\AA)$ consisting of $\KKK$-projective (resp., $\KKK$-injective) $N$-complexes.
A \emph{projective $N$-resolution} (resp., \emph{injective $N$-resolution}) of $X \in\KKK_{N}(\AA)$ is an $N$-quasi-isomorphism $P_X\to X$
(resp., $X\to I_X$) with $P_X\in\KKK_{N}^{\mrm{p}}(\AA) \cap \KKK_{N}(\cat{Prj}\AA)$
(resp., $I_X\in\KKK_{N}^{\mrm{i}}(\AA)\cap \KKK_{N}(\cat{Inj}\AA)$). 
\end{definition}

Clearly $\KKK_{N}^{\mrm{p}}(\AA)$ (resp., $\KKK_{N}^{\mrm{i}}(\AA)$)
is a triangulated subcategory closed under coproducts (resp., products) in 
$\KKK_{N}(\AA)$.
The canonical functor $\KK{\AA} \to \DD{\AA}$ restricts to fully faithful functors 
$\KKK_{N}^{\mrm{p}}(\AA) \to \DD{\AA}$ and $\KKK_{N}^{\mrm{i}}(\AA) \to \DD{\AA}$
by Lemma \ref{fully faithful}. 
By Lemma \ref{lem:localqis1}, 
$\KKK_{N}^{\mrm{-}}(\cat{Prj}\AA)$ (resp., $\KKK_{N}^{\mrm{+}}(\cat{Inj}\AA)$)  is contained in $\KKK_{N}^{\mrm{p}}(\AA)$ (resp., $\KKK_{N}^{\mrm{i}}(\AA)$).

We have the following result which generalizes a classical result for the case $N=2$ \cite{BN,Sp}.

\begin{theorem}\label{cor:subeqv03}
The following hold.
\begin{enumerate}
\item Assume that $\AA$ is an $Ab4$-category with enough projectives.
Then $(\KKK_{N}^{\mrm{p}}(\AA),\KKK_{N}^{\mrm{a}}(\AA))$ is a stable t-structure
in $\KKK_{N}(\AA)$ and we have a triangle equivalence $\KKK_{N}^{\mrm{p}}(\AA) \simeq \DDD_{N}(\AA)$.
Moreover, any object in $\KKK_{N}^{\mrm{p}}(\AA)$ is isomorphic to an object in $\KKK_{N}^{\mrm{p}}(\AA) \cap \KKK_{N}(\cat{Prj}\AA)$, hence every object in $\KKK_{N}(\AA)$ admits a projective $N$-resolution.
\item Assume that $\AA$ is an $Ab4^*$-category with enough injectives.
Then $(\KKK_{N}^{\mrm{a}}(\AA),\KKK_{N}^{\mrm{i}}(\AA))$ is a stable t-structure
in $\KKK_{N}(\AA)$ and we have a triangle equivalence $\KKK_{N}^{\mrm{i}}(\AA) \simeq \DDD_{N}(\AA)$.
Moreover, any object in $\KKK_{N}^{\mrm{i}}(\AA)$ is isomorphic to an object in $\KKK_{N}^{\mrm{i}}(\AA) \cap \KKK_{N}(\cat{Inj}\AA)$, hence every object in $\KKK_{N}(\AA)$ admits an injective $N$-resolution.
\end{enumerate}
\end{theorem}

To prove Theorem \ref{cor:subeqv03}, we need the following easy observation.

\begin{lemma} \label{prop:limithlimit}
Let $\AA$ be an $Ab3$-category, and $f_i: X_{i} \arr X_{i+1}~(i=0,1,\cdots )$ a sequence of
morphisms in $\cat{C}_{N}(\AA)$. Assume that each $j \in \mbb{Z}$ admits some $n \in \mbb{N}$
such that $f^j_i: X^{j}_{i} \arr X^{j}_{i+1}$ is a split monomorphism for $i \geq n$.
Then we have an exact sequence
$0 \arr{\coprod}_{i\ge0}X_i \xarr{1-\coprod_if_i} 
{\coprod}_{i\ge0}X_i \arr \varinjlim X_i \arr 0$
in $(\cat{C}_{N}(\AA),\mcal{S}_{N}(\AA))$ for the inductive limit $\varinjlim X_i$ in $\cat{C}_{N}(\AA)$.
Therefore $\varinjlim X_i$ is isomorphic to the homotopy colimit
$\underset{\lgarr}{\opn{hlim}}\ X_i$ in $\cat{K}_{N}(\AA)$.
\end{lemma}

\begin{proof}
We have a split exact sequence
$0 \arr{\coprod}_{i\ge0}X^{j}_i \xarr{1-\coprod_if^j_i} {\coprod}_{i\ge0}X^{j}_i \arr \varinjlim X^{j}_i \arr 0$
in $\AA$ for any $j$ by our assumption. Thus the assertions follow.
\end{proof}

\begin{proof}[of Theorem \ref{cor:subeqv03}]
We only prove (1) since (2) is the dual.
By Lemma \ref{fully faithful},  it is enough to show 
$\KKK_{N}(\AA) = \KKK_{N}^{\mrm{p}}(\AA) * \KKK_{N}^{\mrm{a}}(\AA)$ to prove the first statement.
\par\noindent
For a complex $X \in \KKK_{N}(\AA)$, we shall construct an $N$-quasi-isomorphism $s: P \to X$ with $P \in \KKK_{N}^{\mrm{p}}(\AA)\cap \KKK_{N}(\cat{Prj}\AA)$. 
Applying Lemma \ref{prop:limithlimit} to a sequence
$\iota_i:{\sigma}_{\leq i}X\to {\sigma}_{\leq i+1}X$ of morphisms, we have $X=\varinjlim X_i
\simeq\underset{\lgarr}{\opn{hlim}}\ X_i$ in $\KKK_N(\AA)$.
By Theorem \ref{cor:subeqv}, there is an $N$-quasi-isomorphism 
$s_i:P_i \arr {\sigma}_{\leq i}X$ with $P_{i} \in \KKK_{N}^{-}(\cat{Prj}\AA)$. 
Since the mapping cone $\opn{C}(s_{i+1})$ is $N$-acyclic, 
by  Lemma \ref{lem:localqis1}
we have a commutative diagram in $\KKK_{N}(\AA)$
\[\xymatrix@R1em{
P_i \ar[r]^{s_i}\ar[d]^{f_i}& {\sigma}_{\leq i}X\ar[d]^{\iota_i}\\
P_{i+1}\ar[r]_{s_{i+1}} & {\sigma}_{\leq i+1}X\ar[r] &\opn{C}(s_{i+1}).
}\]
Therefore we have a morphism between triangles in $\KKK_{N}(\AA)$
\[\xymatrix@C=4em@R1.5em{
{\coprod}_iP_i\ar[r]^{1-\coprod_if_i}\ar[d]^{{\coprod}_is_i}
& {\coprod}_iP_i \ar[r]^{u} \ar[d]^{{\coprod}_is_i}
& P \ar[r]^{v} \ar[d]^{s}& \Sigma{\coprod}_iP_i \ar[d]^{\Sigma{\coprod}_is_i}\\
{\coprod}_i{\sigma}_{\leq i}X\ar[r]_{1-\coprod_i\iota_i} 
& {\coprod}_i{\sigma}_{\leq i}X\ar[r] 
 & X \ar[r] & \Sigma{\coprod}_i{\sigma}_{\leq i}X.
}\]
Since $\AA$ is $Ab4$, ${\coprod}_is_i$ is an $N$-quasi-isomorphism, 
hence so is $s$.
The upper triangle shows $P\in \KKK_{N}^{\mrm{p}}(\AA)\cap \KKK_{N}(\cat{Prj}\AA)$.

Now we prove the second statement.
For any $X\in\KKK_{N}^{\mrm{p}}(\AA)$, the above construction gives a triangle $P\xrightarrow{s} X\to Y\to P[1]$ in $\KKK_N(\AA)$ with $P\in\KKK_{N}^{\mrm{p}}(\AA)\cap \KKK_{N}(\cat{Prj}\AA)$ and $Y\in \KKK_{N}^{\mrm{a}}(\AA)$. Since $\KKK_N^{\mrm{p}}(\AA)$ is a triangulated subcategory of $\KKK_N(\AA)$, we have $Y\in\KKK_N^{\mrm{a}}(\AA)\cap\KKK_N^{\mrm{p}}(\AA)$.
Thus $Y\simeq 0$ and hence $s$ is an isomorphism in $\KKK_N(\AA)$.
\end{proof}

\begin{remark}\label{rem:spcpx}
Later we need a slightly more general version of Theorem \ref{cor:subeqv03} as follows.
\par\noindent
Let $\AA$ be an $Ab4$-category with enough projectives and
$\PP$ an additive subcategory of $\Proj\AA$ closed under coproducts
such that any object in $\Proj\AA$  is an epimorphic image from some object of $\PP$.
Then the proof of Proposition \ref{prop:localqis2} gives triangle equivalences
\[
\KKK_{N}(\PP)\cap\KKK_{N}^{\mrm{p}}(\AA) \simeq\DDD_{N}(\AA)\ \mbox{ and }\ 
\KKK_{N}^{-}(\PP) \simeq\DDD_{N}^{-}(\AA).
\]
For example, the category $\cat{Free} R$ of free modules over a ring $R$ satisfies this condition.
\end{remark}

\begin{example}\label{KW example}
Take a projective $2$-resolution $\cdots \xarr{d^{-2}} P^{-1} \xarr{d^{-1}} P^{0}$ of $X\in\AA$.
Then a projective $N$-resolution of $X$ is given by the following.
\[\xymatrix@C=1.5em@R=0em{
\mbox{degree}: & \scriptstyle{-N-1}&\scriptstyle {-N}& \scriptstyle{-N+1}& \scriptstyle{-N+2}&  & \scriptstyle{-1}& \scriptstyle{0}& \scriptstyle{1}& \scriptstyle{2}\\
{P_X}:\cdots \ar[r] ^{1}& {P^{-3}} \ar[r]^{d^{-3}} & {P^{-2}} \ar[r]^{d^{-2}} & {P^{-1}} \ar[r]^{1} & {P^{-1}} \ar[r]^{1} & \cdots \ar[r]^{1} & P^{-1} \ar[r]^{d^{-1}} & {P^{0}} \ar[r] & 0 \ar[r] & 0 \ar[r]  &\cdots.
} \] 
Although the 2-complex $\cdots \xarr{d^{-2}} P^{-1} \xarr{d^{-1}} P^{0}\xarr{d^0} X \to 0$ is 2-acyclic for some $d^0:P^0\to X$,
the $N$-complex $Y$ below is not $N$-acyclic for $N>2$ since $\opn{H}^{1}_{(1)} (Y) \simeq X$. 
On the other hand, the following $N$-complex $Z$ is $N$-acyclic.
The truncation $\tau_{\le 0}Z$ is not a projective $N$-resolution of $X$, but that of $\Sigma\Theta^{-1}(X)=\mu_{N-1}^{0}(X)$ 
since we have a triangle
$\Theta^{-1}X\to Z\to\tau_{\le 0}Z\to \Sigma\Theta^{-1}X$.
\[\xymatrix@C=1.5em@R=0em{
\mbox{degree}: & \scriptstyle{-N-1}&\scriptstyle {-N}& \scriptstyle{-N+1}& \scriptstyle{-N+2}&  & \scriptstyle{-1}& \scriptstyle{0}& \scriptstyle{1}& \scriptstyle{2}\\
{Y}:\cdots \ar[r] ^{1}& {P^{-3}} \ar[r]^{d^{-3}} & {P^{-2}} \ar[r]^{d^{-2}} & {P^{-1}} \ar[r]^{1} & {P^{-1}} \ar[r]^{1} & \cdots \ar[r]^{1} & P^{-1} \ar[r]^{d^{-1}} & {P^{0}} \ar[r]^{d^0} & {X} \ar[r] & 0\ar[r]  &\cdots \\
{Z}:\cdots \ar[r] ^{1}& {P^{-2}} \ar[r]^{1} & {P^{-2}} \ar[r]^{d^{-2}} & {P^{-1}} \ar[r]^{d^{-1}} &{P^{0}} \ar[r]^{1} & \cdots \ar[r]^{1} & {P^{0}} \ar[r]^{1} & 
{P^{0}} \ar[r]^{d^0} & {X} \ar[r] & 0 \ar[r]  &\cdots } \] 
\end{example}

Let $\MM$ be a full subcategory of $\AA$.
We denote by $\CCC_{N, \MM}( \AA)$ the full subcategory of $\CCC_N(\AA)$ consisting of
$X$ such that $\opn{H}_{(r)}^{i}(X) \in \MM$ for any $0 < r < N$ and $i\in \mathbb{Z}$.
Then $\KKK_{N, \MM}(\AA)$ and $\DDD_{N, \MM}(\AA)$ denote 
the corresponding full subcategories of 
$\KKK_{N}(\AA)$ and $\DDD_{N}(\AA)$ respectively. 
In the case that $\MM$ is a Serre subcategory, that is, closed under subobjects, quotient objects and extensions, 
then $\KKK_{N, \MM}(\AA)$ (resp., $\DDD_{N, \MM}(\AA)$) 
is a thick subcategory of $\KKK_{N}(\AA)$ (resp., $\DDD_{N}(\AA)$). 
We use the notations
$\CCC_{N, \MM}^{ \sharp , \natural}( \AA) = \CCC_{N}^{\sharp , \natural}( \AA) \cap\CCC_{N, \MM}( \AA)$, 
$\KKK_{N, \MM}^{ \sharp , \natural}( \AA) = \KKK_{N}^{\sharp , \natural}( \AA) \cap\KKK_{N, \MM}( \AA)$  and 
$\DDD_{N, \MM}^{\sharp, \natural}(\AA)=\DDD_{N}^{\sharp,\natural}(\AA)\cap\DDD_{N, \MM}(\AA)$ 
for $\sharp =$nothing$,-,+,\bo$ and $\natural=$nothing$,-,+,\bo$. 
By Proposition \ref{b+-}, we have $\DDD_{N, \MM}^{\sharp, \bo}(\AA)\simeq\DDD_{N, \MM}^{\bo}(\AA)$ etc.

\begin{proposition} \label{prop:localqis3}
Let $\MM$ be an additive full subcategory of $\AA$ satisfying ($V_1$) and ($V_2$).
\begin{enumerate}
\item
For any $X \in \CCC_{N, \MM}^{-}( \AA)$, there is an $N$-quasi-isomorphism $L \to X$ with $L \in \CCC_{N}^{-}(\MM)$.
\item
$\KKK_{N, \MM}^{-,\natural}( \AA)   \subset \KKK_N^{-,\natural} (\MM )* \KKK_N^{-,\mrm{a}}(\AA)$ for $\natural=$nothing$,\bo$.
\end{enumerate}
\end{proposition}

\begin{proof}
(1) There exists $n_0$ such that $X^i =0$ for any $i > n_0$. 
{Set} $L_{n_0}=X$.
We shall construct a sequence of $N$-quasi-isomorphisms 
$v_n : L_{n-1} \to L_{n}$ in  $\CCC_N ^{}(\AA)$ for $n\le n_0$ such that  
\[
L_n ^i \in \MM \quad (i > n),\ \opn{B}^i _{(r)} (L_n ) \in \widehat{\MM} \quad (i>n,\ 0<r<N)\ \mbox{ and }v_n^i=\opn{id} \quad (i>n)\] 
Then we get an $N$-quasi-isomorphism $L=\underset{\lglarr}{\lim}\,L_{n}\to X$ with $L\in\CCC^{-}_N(\MM)$. 
Suppose $n < n_0$ and let $L_n$ satisfy the conditions above. 
The exact sequence
$0 \to \opn{H}^{n}_{(1)} (L_n ) \to \opn{C}^{n} _{(N-1)} (L_n) \to \opn{B}^{n+1} _{(1)} (L_n) \to 0$
implies $\opn{C}^{n} _{(N-1)} (L_n) \in \widehat{\MM}$.  
Applying Lemma \ref{set01}(2) to the canonical epimorphism 
$\rho : L^{n} _n \to \opn{C}^{n} _{(N-1)} (L_n)$, 
we get a morphism 
$v: M \to L_n ^{n}$ with $M \in \MM$ such that 
$\rho v$ is an epimorphism. 
Set $L_{n-1} = \opn{V}_n(L_n,v)$ and $v_{n-1}= p_n(v)$. 
{\[ \xymatrix@C=4em@R1.5em{
\ar @{} [dr] |{\qquad (E)} L_{n-1}^{n-N+1}  \ar[r] ^{d^{\{N-1\}}_{L_{n-1}}} \ar[d]^{v^{n-N+1}_n} 
& M=L_{n-1}^{n}  \ar[d]^{v}  \\
L_{n}^{n-N+1}  \ar[r]_{d^{\{N-1\}}_{L_{n}}}
& L_{n}^{n}  \ar[r]_{\rho} 
&  \opn{C}^{n} _{(N-1)} (L_n)} \] }
Since $\rho$ is the cokernel of $d^{\{N-1\}}_{L_{n}}$ and $\rho v$ is an epimorphism,
$(v\ d^{\{N-1\}}_{L_{n}}):L_{n-1}^{n}\oplus L_{n}^{n-N+1} \to L_n^n$ is an epimorphism, 
which shows $(E)$ is an exact square. 
Thus $v_{n-1} = p_n(v)$ is an $N$-quasi-isomorphism 
by Definition-Proposition \ref{elmap02}.

Now we show that $\opn{B}^i _{(r)} (L_{n-1} ) \in \widehat{\MM}$ for any $i>n-1$ and $0<r<N$. 
If $i>n$, then $\opn{H}^{i} _{(N-r)} (L_{n-1}) = \opn{H}^{i} _{(N-r)} (L_{n}) \in \MM$ holds.
Moreover $\opn{Z}^{i} _{(N-r)} (L_{n-1}) = \opn{Z}^{i} _{(N-r)} (L_{n})$ belongs to $\widehat{\MM}$
since $0\to\opn{Z}^{i} _{(N-r)} (L_{n})\to L^{i}_{n}\to\opn{B}^{i+N-r}_{{(N-r)}}(L_{n})\to0$ is exact.
Therefore $\opn{B}^{i}_{(r)} (L_{n-1} )  \in \widehat{\MM}$ holds.
To see $\opn{B}^n_{(r)}  (L_{n-1} )\in \widehat{\MM}$, it suffices to show $\opn{C}^{n}_{(r)} (L_{n-1} )  \in \widehat{\MM}$
since $L_{n-1}^n \in \MM$. 
But this is clear since  $\opn{B}^{n+N-r} _{(N-r)} (L_{n-1}) =\opn{B}^{n+N-r} _{(N-r)} (L_{n}) \in \widehat{\MM}$ and 
$\opn{H}^{n}_{(N-r)} (L_{n-1} ) \in \MM$. 
\par\noindent
(2) For given $X\in \KKK_{N, \MM}^{-}( \AA)$, there is an $N$-quasi-isomorphism
$s:L \to X$ with $L \in \KKK_{N} ^{-} (\MM)$ by (1). 
We get the first inclusion since $\mrm{C}(s) \in \KKK_{N}^{\mrm{a}}(\AA)$. 
If $X \in \KKK_{N, \MM}^{-}( \AA)$, the construction shows $\mrm{C}(s) \in\KKK_{N}^{-,\mrm{a}}(\AA)$. 
If $X \in \KKK_{N, \MM}^{-, \bo}( \AA)$, then obviously we have $L \in \KKK_{N} ^{-, \bo} (\MM)$. 
\end{proof}

\begin{theorem}\label{Serre}
If $\MM$ is a Serre subcategory satisfying the condition ($V_1$), then $\DDD_{N}^{\natural}(\MM)\simeq\DDD_{N, \MM}^{\natural}(\AA)$ for $\natural=\bo, -$.
\end{theorem}

\begin{proof} 
Since $\MM$ is a Serre subcategory, it satisfies the condition ($V_2$) and we have $ \KKK_N^{-, \natural} (\MM ) \subset \KKK^{-, \natural} _{N, \MM} (\AA)$. 
By Proposition \ref{prop:localqis3}(2), we have $\KKK^{-, \natural} _{N, \MM} (\AA ) = \KKK_N^{-, \natural } (\MM )* \KKK_N^{-,\mrm{a}}(\AA)$.
Applying Lemma \ref{fully faithful} to $\UU=\KKK^{-,\natural}_N(\MM)$ and $\VV=\KKK^{-,\mrm{a}}_N(\AA)$, we have triangle equivalences
$\DDD^{\natural}_N(\MM)\simeq\frac{\KKK^{-,\natural}_N(\MM)}{\KKK^{-,\mrm{a}}_N(\MM)}
=\frac{\UU}{\UU\cap\VV}\simeq\frac{\UU*\VV}{\VV}=
\frac{\KKK^{-,\natural}_{N,\MM}(\AA)}{\KKK^{-,\mrm{a}}_N(\AA)}\simeq\DDD^{\natural}_{N,\MM}(\AA)$ as desired.
\end{proof}

\subsection{Homotopy categories of injective objects}\label{HCIO}

In this subsection, we shall show that $\KKK_{N}(\Ij \AA)$ is compactly generated if $\AA$ satisfies some conditions. 

An  \emph{$Ab5$-category} is an $Ab3$-category that has exact filtered colimits. A \emph{Grothendieck category} is an $Ab5$-category with a generator.
A Grothendieck category $\AA$ is called \emph{locally noetherian} if $\AA$ has a generating set of noetherian objects.
In this case, $\Inj\AA$ is closed under arbitrary coproducts \cite[Theorem 8.7]{Po}, and therefore the triangulated category $\KKK_N(\Inj\AA)$ has arbitrary coproducts.

For an additive category $\BB$ with arbitrary coproducts, an object $C$ is called {\it compact} in $\BB$
if the canonical morphism $\coprod_i\Hom_{\BB}(C, X_i)\iso \Hom_{\BB}(C, \coprod_iX_i)$ is an isomorphism for any coproduct $\coprod_iX_i$ in $\BB$.
We denote by $\BB^{\co}$ the category of compact objects in $\BB$. 
A triangulated category $\mcal{D}$ with arbitrary coproducts
is called \emph{compactly generated} by a set $\cat{S}$ of compact objects if any non-zero object of $\mcal{D}$ has a non-zero morphism from a shift of some object of $\cat{S}$.

Let $\nt\AA$ be the subcategory of $\AA$ consisting of noetherian objects.
For a locally noetherian Grothendieck category $\AA$, it is easy to see $\nt\AA$ is a skeletally small 
Serre subcategory satisfying ($V_1$) and ($V_2$).
By Theorem \ref{Serre}, we can identify
$\DDD_{N}^{\bo}(\cat{noeth} \AA)$ with $\DDD_{N, {\cat{noeth}\AA}}^{\bo}(\AA)$.

We aim to prove the $N$-complex version of a result of Krause \cite{Kr2}.

\begin{theorem}\label{Kijcp2}
Let $\AA$ be a locally noetherian Grothendieck category.
Then $\KKK_{N}(\Ij \AA)$ is a compactly generated triangulated category such that 
the canonical functor $\KKK_{N}(\Ij \AA) \to \DDD_{N}(\AA)$ induces
an equivalence between $\KKK_{N}(\Ij \AA)^{\co}$ and $\DDD_{N}^{\bo}(\cat{noeth} \AA)$.
\end{theorem}

In the rest, $\AA$ is a locally noetherian Grothendieck category. 
Recall that $I_X\in\KKK_{N}^{\mrm{i}}(\cat{Inj}\AA)$ stands for the injective $N$-resolution of an object $X$ in 
$\KKK_{N}(\AA)$.

\begin{lemma}(cf. \cite[Lemma 2.1]{Kr2}) \label{Kijcp}
The object $I_{\mu_{r}^{s}(M)}$ is compact in $\KKK_N(\Ij \AA)$ for any $M \in \nt\AA$, $s\in\z$ and $0<r<N$.
\end{lemma}

\begin{proof}
For any $Y \in \KKK_N(\Ij\AA)$, we have the following isomorphisms for sufficiently small $t$:
{\small\[\begin{aligned}
\Hom_{\KKK_{N}(\AA)}(I_{\mu_{r}^{s}(M)}, Y) & \simeq \Hom_{\KKK_{N}(\AA)}(I_{\mu_{r}^{s}(M)},  \tau_{\geq t}Y) & \simeq \Hom_{\KKK_{N}(\AA)}({\mu_{r}^{s}(M)}, \tau_{\geq t}Y) \\ & \simeq \Hom_{\KKK_{N}(\AA)}({\mu_{r}^{s}(M)}, Y). 
\end{aligned}\]}
The first and third isomorphisms come from  $I_{\mu_{r}^{s}(M)}, ~{\mu_{r}^{s}(M)} \in \KKK_N^+(\AA)$ and the second one from Lemma \ref{lem:localqis1}.
Also we have 
$\Hom_{\KKK_{N}(\AA)}(\mu_{r}^{s}(M),Y)\simeq\opn{H}^{s-r+1}_{(r)}(\Hom_{\AA}(M,Y))$
by \eqref{homfrommu1}. This completes the proof since $M\in\nt\AA$ is compact in $\AA$.
\end{proof}

Let 
$\cat{S}$ stand for a set of representatives of isomorphism classes of objects $\{I_{\mu_r^{s}(M)}\mid M\in\nt\AA,\ s\in\z, 0<r <N-1\}$
in $\KKK_{N}(\Ij \AA)$.

\begin{lemma}\label{lem:Kr2.2}(cf. \cite[Lemma 2.2]{Kr2})
$\KKK_{N}(\Ij \AA)$ is compactly generated by $\cat{S}$.
\end{lemma}

\begin{proof}
By Lemma \ref{Kijcp}, any object of $\cat{S}$ is compact in $\KKK_{N}(\Ij \AA)$. Let $X \in \KKK_{N}(\Ij \AA)$ be a non-zero object.
Assume that $\opn{H}^{i}_{(r)}(X)\not=0$ for some $i\in\z$ and $0<r<N$.
Since $\AA$ is locally noetherian, there is a non-zero morphism 
$M \to \opn{Z}^{i}_{(r)}(X) \to \opn{H}^{i}_{(r)}(X)$ with $M\in\nt\AA$.
Using the commutative diagram in Lemma \ref{homfrommu}(1), we have $\Hom_{\KKK_{N}(\AA)}(\mu_{r}^{i+r-1}(M),X)\not=0$.
\par\noindent
Assume that $X$ is $N$-acyclic.
Since $X\not=0$ in $\KKK_{N}(\Ij \AA)$, there are 
$i\in\z$ and $0<r<N$ with $\opn{Z}^{i}_{(r)}(X) \not\in \Ij \AA$ by Lemma \ref{0cpx}(3). 
Baer criterion \cite[Lemma A10]{Kr3} gives an object $M$ of $\cat{noeth}\AA$ with 
$\Ext_{\AA}^{1}(M,\opn{Z}^{i}_{(r)}(X))\not=0$, 
which implies $\Hom_{\KKK_{N}(\AA)}(\mu^{i+N-1}_{N-r}(M),X)\not=0$ by Lemma \ref{homfrommu}(3).
\end{proof}

Now we are ready to prove Theorem \ref{Kijcp2}.

\begin{proof}[of Theorem \ref{Kijcp2}]
Lemma  \ref{lem:Kr2.2} implies $\KKK_{N}(\Ij \AA)=\Loc \cat{S}$  (see \cite[1.6]{Ne1}). 
Hence by \cite[Lemma 2.2]{Ne1}, $\KKK_{N}(\Ij \AA)^{\co}$ coincides with $\thick\cat{S}$.
On the other hand, the equivalence $\KKK_{N}^{\mrm{i}}(\cat{Inj}\AA) \simeq \DDD_{N}(\AA)$ in Theorem \ref{cor:subeqv}(2) yields 
$\thick_{\KKK_N^{\mrm{i}}(\Inj \AA)}\cat{S}\simeq \\ \thick_{\DDD_N(\AA)}(\cat{noeth} \AA)\simeq \DDD_N^{\bo}(\cat{noeth} \AA)$.
\end{proof}

\subsection{Derived functor}\label{DFunNcpx}

In this subsection, we study the derived functor of a triangle functor $\KKK_{N}(\AA) \to
\KKK_{N'}(\AA')$ for abelian categories $\AA$, $\AA'$.

\begin{definition}\label{derf}
Let $\TT$ be a triangulated category, $\UU$ a full triangulated subcategory of $\TT$ and $Q:\TT\to\TT/\UU$ the canonical functor.
For a triangle functor $F:\TT \to \TT'$, the \emph{right derived functor} (resp., \emph{left derived functor}) of $F$
with respect to $\UU$ is a triangle functor
\[\begin{aligned}
\bsym{R}_{\UU}F:\TT/\UU \to \TT' &\quad
(\text{resp.,}\ \bsym{L}_{\UU}F:\TT/\UU \to \TT'\  )
\end{aligned}\]
together with a functorial morphism of triangle functors
\[\begin{aligned}
\xi :F \to (\bsym{R}_{\UU}F) Q & \quad
(\text{resp.,}\ \xi :(\bsym{L}_{\UU}F) Q \to F \  )
\end{aligned}\]
with the following property:\par\noindent
For a triangle functor $G:\TT/\UU \to \TT'$ and
a functorial morphism of triangle functors
$\zeta:F\to G Q$
(resp., $\zeta:G Q \to F$),
there exists a unique functorial morphism $\eta:\bsym{R}_{\UU}F\to G$ 
(resp., $\eta: G \to \bsym{L}_{\UU}F$)
of triangle functors such that
$\zeta=(\eta Q)\xi$ (resp.,\ $\zeta=\xi (\eta Q)$).
\[\begin{xy}
(0,0)*+{\TT}="a1",
(0,-14)*+{\TT/\UU}="a2",
(30,0)*+{\TT'}="a3",
(15,-7)*+{ }="a5",
(19,-15)*+{ }="a6",
\ar@{->}_Q "a1";"a2",
\ar@{->}^F "a1";"a3",
\ar@{->}^{\bsym{R}_{\UU}F} "a2";"a3",
\ar@/_25pt/@{->}_G "a2";"a3",
\ar@{=>} "a5";"a6"
\end{xy}\]
\end{definition}

We recover a classical Existence Theorem of derived functors as follows:

\begin{theorem}[Existence Theorem] \label{thm:exderfun}
Let $\TT$ be a triangulated category, $\UU$ its full triangulated subcategory, and $Q:\TT\to\TT/\UU$ the canonical functor.
For a triangle functor $F:\TT \to \TT'$, assume that there exists a full triangulated subcategory $\VV$ of $\TT$ such that $\TT=\UU*\VV$
and $F(\UU\cap\VV) = \{0\}$.
Then there exists the right derived functor $(\bsym{R}_{\UU}F, \xi)$ of $F$ with respect to $\UU$ such that
$\xi_X: FX \arr (\bsym{R}_{\UU}F)QX$  is an isomorphism for $X\in\VV$.
\end{theorem}

\begin{proof}
Let $Q_1:\TT\xrightarrow{}\TT/(\UU\cap\VV)$ and $Q_2:\TT/(\UU\cap\VV)\to\TT/\UU$ be the canonical functors. Then $Q=Q_2Q_1$ holds.
Since $F(\UU\cap\VV)=0$, the functor $F:\TT\to\TT'$ factors as $\TT\xrightarrow{Q_1}\TT/(\UU\cap\VV)\xrightarrow{F'}\TT'$ by universality.
By Lemma \ref{fully faithful}, the functor $Q_2:\TT/(\UU\cap\VV)\to\TT/\UU$ has a right adjoint $R:\TT/\UU\to\TT/(\UU\cap\VV)$.

We shall show that $\bsym{R}_{\UU}F=F'R$ satisfies the condition. 
We have only to give a functorial isomorphism $\Hom_\triangle(F,GQ)\simeq\Hom_\triangle(F' R ,G)$ for any triangle functor $G:\TT/\UU\to\TT'$, where $\Hom_\triangle$ is the class of morphisms between triangle functors.
Indeed, we have $\Hom_\triangle(F,GQ)\simeq\Hom_\triangle(F',GQ_2)$ by \cite[Proposition 3.4]{Har}, 
and $\Hom_\triangle(F',GQ_2)\simeq\Hom_\triangle(F'R,G)$ by a triangle functor version of \cite[Proposition X.7.3]{Mc}. 
\end{proof}

We apply these to the setting of $N$-complexes.

\begin{definition}[Derived Functor]
Let $\AA$ and $\mcal{A'}$ be abelian categories, 
and $F: \KKK_{N}^{\natural}(\AA) \arr \KKK_{N'}(\mcal{A'})$ a triangle functor
where $\natural=$nothing$,-,+,\bo$.
We define the \emph{right} (resp., \emph{left}) \emph{derived functor} of $F$ as
\[\begin{aligned}
\bsym{R}^{\natural}F=\bsym{R}_{\UU}(Q'F):\DDD_N^{\natural}(\AA) \to \DDD_N(\AA') &\quad
(\text{resp.,}\ \bsym{L}^{\natural}F=\bsym{L}_{\UU}(Q'F):\DDD_N^{\natural}(\AA) \to \DDD_N(\AA')),
\end{aligned}\]
where $Q':\KKK_N(\AA')\to\DDD_N(\AA')$ is the canonical functor,
$\TT=\KKK_N^{\natural}(\AA)$ and $\UU=\KKK_N^{\natural, \mrm{a}}(\AA)$.
\end{definition}

According to Theorems \ref{cor:subeqv}, \ref{cor:subeqv03} and \ref{thm:exderfun},
we have the following $N$-complex version of classical results \cite{Har,BN,Sp}.

\begin{corollary}\label{thm:exderfun02}
Let $\AA$ and $\mcal{A'}$ be abelian categories, and $F : \KKK_{N}(\AA) \arr \KKK_{N'}(\mcal{A'})$ a triangle functor. 
Then the following hold.
\begin{enumerate}
\item If $\AA$ has enough injectives, then $\bsym{R}^{+}F:\DDD_{N}^{+}(\AA) \arr \DDD_{N'}(\mcal{A'})$ exists.
\item If $\AA$ has enough projectives, then $\bsym{L}^{-}F:\DDD_{N}^{-}(\AA) \arr \DDD_{N'}(\mcal{A'})$ exists.
\item If $\AA$ is an $Ab4^*$-category with enough injectives, then
$\bsym{R}F:\DDD_{N}(\AA) \arr \DDD_{N'}(\mcal{A'})$ exists.
\item If $\AA$ is an $Ab4$-category with enough projectives, then
$\bsym{L}F:\DDD_{N}(\AA) \arr \DDD_{N'}(\mcal{A'})$ exists.
\end{enumerate}
\end{corollary}

We end this subsection with considering $\opn{Ext}$ and $\opn{Tor}$ groups. 
As we will see in Proposition \ref{ExtTor_comparison}, these homology groups are related to classical $\opn{Tor}$ and $\opn{Ext}$. 

\begin{definition}\label{ExtTor}
Let $A$ be a ring, $X$ a right $A$-module and $Y$ a left $A$-module.
We have triangle functors $\Hom_A(X,-):\KKK_N(\Mod A)\to\KKK_N(\Mod\z)$ and $-\otimes_AY:\KKK_N(\Mod A)\to\KKK_N(\Mod\z)$.
By Corollary \ref{thm:exderfun02}, we have derived functors
\[\bsym{R}\Hom_A(X,-):\DDD_N(\Mod A)\to\DDD_N(\Mod\z)\ \mbox{ and }\ -\lten_AY:\DDD_N(\Mod A)\to\DDD_N(\Mod\z).\] 

For a right $A$-module $Z$, $n\in\z$ and $0<r<N$, set
\[{}_{r}\opn{Ext}_{A}^n(X,Z)=\opn{H}_{(r)}^{n}(\bsym{R}\Hom_A(X,Z))\ \mbox{ and }\ 
{}_{r}\opn{Tor}_n^{A}(Z,Y)=\opn{H}^{-n}_{(r)}(Z\lten_AY).\]
\end{definition}

\begin{proposition}\label{ExtTor_comparison}
We have the following isomorphisms for $i\geq 0$ and $0<r<N$.
\begin{enumerate}
\item ${}_{r}\opn{Tor}_{iN}^{A}(X,Y)=\opn{Tor}_{2i}^{A}(X,Y)$ and ${}_{r}\opn{Ext}^{iN}_{A}(X,Z)=\opn{Ext}^{2i}_{A}(X,Z)$.
\item ${}_{r}\opn{Tor}_{iN+s}^{A}(X,Y)=
\left\{\begin{array}{ll}
\opn{Tor}_{2i+1}^{A}(X,Y)&r=s.\\
0 &r\neq s
\end{array}\right.$
\item ${}_{r}\opn{Ext}^{iN+s}_{A}(X,Z)=
\left\{\begin{array}{ll}
\opn{Ext}^{2i+1}_{A}(X,Z)&r=N-s.\\
0&r\neq N-s
\end{array}\right.$
\end{enumerate}
\end{proposition}

\begin{proof}
We give a proof only for $\opn{Tor}$.
Let $\cdots \xarr{d^{-2}} P^{-1} \xarr{d^{-1}} P^{0} \xarr{d^0} Y \to 0$ be a projective $2$-resolution of $Y\in\Mod A^{\rm op}$.
We have a projective $N$-resolution of $Y$ by Example \ref{KW example}:
\[
\stackrel{\mbox{degree}}{\cdots} \to \stackrel{-N-2}{P^{-3}} \xrightarrow{1} 
\stackrel{-N-1}{P^{-3}} \xrightarrow{d^{-3}} \stackrel{-N}{P^{-2}} \xrightarrow{d^{-2}} 
\stackrel{-N+1}{P^{-1}} \xrightarrow{1} \cdots \xrightarrow{1} 
 \stackrel{-1}{P^{-1}} \xrightarrow{d^{-1}} \stackrel{0}{P^{0}}.
 \]
Applying $X\otimes_A-$, we can justify the assertions.
\end{proof}

Our Definition \ref{ExtTor} is slightly different from Ext and Tor groups introduced by Kassel and Wambst \cite{KW}.
As we discussed in Example \ref{KW example}, their definitions are interpreted as
\[{}_r\opn{Ext}_{A}^n(X,Z)^{\rm KW}=\opn{H}_{(r)}^{n}(\Hom_A(P_{\Sigma\Theta^{-1}X},Z))\ \mbox{ and }\ 
{}_{r}\opn{Tor}_n^{A}(X,Y)^{\rm KW}=\opn{H}^{-n}_{(r)}(P_{\Sigma\Theta^{-1}X}\ten_AY). \]

\section{Triangle equivalence between derived categories}\label{TrieqDN}

In this section, we show that
the derived category $\DDD_{N}(\AA)$ of $N$-complexes is triangle equivalent to
the ordinary derived category $\DDD(\Morph_{N-2}(\AA))$ where 
$\Morph_{N-2}(\AA)$ is the category of sequences of $N-2$ morphisms in $\AA$. 

\begin{definition}\label{smcat}
Let $\BB$ be an additive category.
The category $\Morph_{N-2}(\BB)$ (resp., $\Morph^{\sm}_{N-2}(\BB)$, $\Morph^{\se}_{N-2}(\BB)$) 
is defined as follows.
\begin{enumerate}
\item An object is a sequence of $N-2$ morphisms (resp., split monomorphisms, split epimorphisms) 
$X: X^{1} \xarr{\alpha_X^1}X^{2} \xarr{\alpha_X^2} \cdots \xarr{\alpha_X^{N-2} } X^{N-1}$ in $\BB$. 
\item A morphism from $X$ to $Y$ is an $(N-1)$-tuple $f=(f^1,\cdots , f^{N-1})$ of morphisms $f^i:X^{i} \to Y^{i}$
which makes the following diagram commutative.
\[\xymatrix@R1em{
X^1\ar[r]^{\alpha_X^1}\ar[d]^{f^1}&X^2\ar[r]^{\alpha_X^2}\ar[d]^{f^2}&\cdots\ar[r]^{\alpha_X^{N-3}}&X^{N-2}\ar[r]^{\alpha_X^{N-2}}\ar[d]^{f^{N-2}} &X^{N-1}\ar[d]^{f^{N-1}}\\
Y^1\ar[r]_{\alpha_Y^1}&Y^2\ar[r]_{\alpha_Y^2}&\cdots\ar[r]_{\alpha_Y^{N-3}} &Y^{N-2}\ar[r]_{\alpha_Y^{N-2}} &Y^{N-1}
}\]
\end{enumerate}
\end{definition}

We can identify $\Morph_{N-2}(\BB)$ with a full subcategory of $\CCC_N(\BB)$ (and $\KKK_N(\BB)$) consisting of $N$-complexes concentrated in degrees $1,\ldots,N-1$.
Indeed, we have isomorphisms
\[\Hom_{\Morph^{\sm}_{N-2}(\BB)}(X,Y) =\Hom_{\CCC_{N}(\BB)}(X,Y) =\Hom_{\KKK_{N}(\BB)}(X,Y)\]
for any $X, Y \in \Morph^{\sm}_{N-2}(\BB)$.
As usual, a set $\cat{S}$ of objects in an abelian category $\AA$ is a \emph{set of generators} if 
any object $X\in\AA$ admits an epimorphism from a coproduct of objects in $\cat{S}$ to $X$.

\begin{theorem}\label{DND}
Let $\AA$ be an $Ab3$-category with a small full subcategory {$\CC$}
 of compact projective generators.
Then we have a triangle equivalence
\[
\DDD_{N}(\AA) \simeq\DDD(\Morph_{N-2}(\AA))
\]
which restricts to the identity functor on $\Morph^{\sm}_{N-2}(\CC)$.
\end{theorem}

We start with the following basic observations.

\begin{lemma}\label{find tilting}  
Let $\BB$ be an additive category.
\begin{enumerate}
\item
Assume that $\BB$ is \emph{idempotent complete}, that is, for any $X \in \BB$ and any idempotent $e\in\End_{\BB}(X)$,
there are an object $Y\in\BB$, and morphisms $p: X \to Y$ and $q:Y \to X$ such that $e=qp$ and $pq=1_Y$.
Then for every object $P$ of $\Morph^{\sm}_{N-2}(\BB)$,
there are objects $C_1, \cdots ,C_{N-1}$ of $\BB$ such that $P \simeq \coprod_{i=1}^{N-1}\mu^{N-1}_{i}(C_{i})$.
\item  For any $P,Q  \in \Morph^{\sm}_{N-2}(\BB)$, we have
$\Hom_{\KKK_{N}(\BB)}(P,\Sigma^{j}Q) = 0 \ (j\not=0)$.
\item $\KKK_N^{\bo}(\BB)=\tri\Morph^{\sm}_{N-2}(\BB)$.
\item Assume that $\BB$ has arbitrary coproducts. Then every object in $\Morph^{\sm}_{N-2}(\BB^{\co})$ is compact in $\CCC_N(\BB)$ (resp., $\KKK_N(\BB)$).
\end{enumerate}
\end{lemma}

\begin{proof}
(1) This is clear.\par\noindent
(2) Let $\widetilde{\BB}$ be the idempotent completion of $\BB$ (e.g. \cite[Definition 1.2]{BS}).
Since $\KKK_N(\BB)$ is a full triangulated subcategory of $\KKK_N(\widetilde{\BB})$, 
we can assume that $\BB$ is idempotent complete.
By (1), we have only to consider the case $P=\mu_r^{N-1}(C)$ and $Q=\mu_{r'}^{N-1}(C')$ for $C, C' \in \BB$ and $0<r,r' <N$.
For the case $j=1$, we have $\Sigma\mu_{r'}^{N-1}(C')=\mu_{N-r'}^{N-r'}(C')$ by Lemma \ref{Sigma of mu}(1),
and it is easy to check that any morphism from $\mu_{r}^{N-1}(C)$ to $\mu_{N-r'}^{N-r'}(C')$ is null-homotopic.
Now we consider the case $j\neq0,1$. 
Since $\Sigma^2=\Theta^N$, there is no degree in which both
$\mu_{r}^{N-1}(C)$ and $\Sigma^{j}\mu_{r'}^{N-1}(C')$ have non-zero terms.
Thus we have $\Hom_{\CCC_N(\BB)}(\mu_{r}^{N-1}(C),\Sigma^{j}\mu_{r'}^{N-1}(C'))=0$.
\par\noindent
(3) For any $C\in\BB$ and $0<r<N$, we have a triangle $\mu^r_1(C)\to\mu^{N-1}_{N-r}(C)\to\mu^{N-1}_{N-r-1}(C)\to\Sigma\mu^r_1(C)$ in $\KKK_N^{\bo}(\BB)$.
Thus $\mu^r_1(C)\in\tri\Morph^{\sm}_{N-2}(\BB)$ holds. By Lemma \ref{Sigma of mu}(2), the assertion follows.
\par\noindent
(4) Taking idempotent completion of $\BB$,
it suffices to show that $\mu_r^{N-1}(C)$ is compact in $\CCC_N(\BB)$ (resp. $\KKK_N(\BB)$) for $C\in\BB^{\co}$.
This follows from \eqref{homfrommu1}.
\end{proof}

\begin{definition}\label{defn:tilting}
Let $\TT$ be a triangulated category with arbitrary coproducts.
A small full subcategory $\SS$ of $\TT^{\co}$ is called a \emph{tilting subcategory} if the following conditions are satisfied.
\begin{enumerate}
\item $\Hom_{\TT}(\SS,\Sigma^i\SS)=0$ for any $i\neq0$.
\item If $X\in\TT$ satisfies $\Hom_{\TT}(\SS,\Sigma^iX)=0$ for any $i\in\z$, then $X=0$.
\end{enumerate}
\end{definition}

The following general result by Keller is basic, where we always regard $\SS$ as a full subcategory of $\Mod \SS$ and $\DDD(\Mod \SS)$ by Yoneda embedding. 

\begin{proposition}\label{tilting theorem}
Let $\TT$ be an algebraic triangulated category with arbitrary coproducts and $\SS$ a tilting subcategory. 
Then we have a triangle equivalence $F:\TT\simeq\DDD(\Mod\SS)$, which restricts to the identity functor on $\SS$.
\end{proposition}  

\begin{proof} 
Although this is well-known, we include a proof for convenience of the reader, because of the lack of proper reference in this setting (cf. \cite[Theorem 8.3.3]{Ke3} for the one-object version).
Replacing objects in $\TT$ with their complete resolutions in the Frobenius category (cf. \cite[Theorem 4.3]{Ke1}, \cite[Theorem 7.5]{Kr}),
we obtain a DG category $\RR$ and a triangle functor $G:\TT\to\DDD(\RR)$ satisfying the following conditions.
\begin{itemize}
\item $\opn{H}^0(\RR)=\SS$ and $\opn{H}^i(\RR)=0$ for any $i\neq0$.
\item $G$ commutes with arbitrary coproducts and induces an equivalence
$\SS\to\widehat{\RR}$, where $\widehat{\RR}$ is the full subcategory of $\DDD(\RR)$ consisting of representable DG functors.
\end{itemize}
Then $G$ induces a triangle equivalence $\Loc\SS\to\Loc\widehat{\RR}$.
Since $\Loc\SS=\TT$ and $\Loc\widehat{\RR}=\DDD(\RR)$ hold by Brown representability, $G:\TT\to\DDD(\RR)$ is a triangle equivalence.
\par\noindent
On the other hand, 
DG functors $\sigma_{\le0}(\RR)\to\RR$ and $\sigma_{\le0}(\RR)\to\opn{H}^0(\RR)=\SS$ are quasi-equivalences \cite{Ke2}
where $\sigma_{\le0}(\RR)$ is the DG category with the same objects as $\RR$ and the morphism spaces given as 
$\Hom_{\sigma_{\le0}(\RR)} (X,Y) = \sigma_{\le 0} \Hom _{\RR} (X,Y)$. 
Hence we have triangle equivalences $\DDD(\RR)\simeq\DDD(\sigma_{\le0}(\RR))\simeq\DDD(\Mod\SS)$ by 
{\cite[9.1]{Ke1}} (cf. \cite[Lemma 3.10]{Ke2}). Thus the assertion follows.
\end{proof}

We need the following general observation.

\begin{proposition}\label{Aus2}
Let $\AA$ be an $Ab3$-category with a small full subcategory $\CC$ of compact projective generators.
Then we have an equivalence $\AA\simeq\Mod\CC$ given by $X\mapsto\Hom_{\AA}(-,X)|_{\CC}$.
In particular, $\mcal{A}$ is a Grothendieck category which satisfies the condition $Ab4^*$.
\end{proposition}

\begin{proof}
See \cite[Chapter IV, Theorem 5.3]{Mit} and \cite[3.4]{Po}. 
\end{proof}

Now we give the following crucial results.

\begin{proposition}\label{find tilting2}
Let $\AA$ be an $Ab3$-category with a small full subcategory 
$\CC$ of compact projective generators.
\begin{enumerate}
\item $\DDD_N(\AA)$ has a tilting subcategory $\Morph^{\sm}_{N-2}(\CC)$. 
\item We have a triangle equivalence $\DDD_N(\AA)\simeq\DDD(\Mod(\Morph^{\sm}_{N-2}(\CC)))$, which restricts to the identity functor on $\Morph^{\sm}_{N-2}(\CC)$.
\end{enumerate}
\end{proposition}

\begin{proof}
(1) 
Set $\SS= \Morph^{\sm}_{N-2}(\CC)$. 
Lemma \ref{find tilting}(4) gives $\SS \subset\KKK_N^{\mrm{p}}(\Proj\AA)^{\co}\simeq\DDD_{N}(\AA)^{\co}$.
Also, $\SS$ satisfies (1) of Definition \ref{defn:tilting} by Lemma \ref{find tilting}(2).  
To show (2) of Definition \ref{defn:tilting}, let $X$ be a non-zero object in $ \DDD _{N} (\AA)$. 
It suffices to find some $C\in \CC$ and $r,s \in \z$ with $0<r<N$ 
such that $\Hom_{\DDD(\AA)}(\mu^{s}_r(C),X)\neq0$. 
Indeed, there exist $i\in\z$ and $0<r<N$ such that $\opn{H}^i_{(r)}(X)\neq0$. 
Since $\CC$ generates $\AA$, we have $\Hom_{\AA}(C,\opn{H}^i_{(r)}(X))\neq0$ 
for some $C \in \CC$. 
So $\Hom_{\DDD(\AA)}(\mu^{i+r-1}_r(C),X)=
\Hom _{\KKK_{N} (\AA)} ( \mu ^{i+r-1}_r (C), X ) 
\neq 0$ by Lemma \ref{homfrommu}(2). 
\par\noindent
(2) This is immediate from (1) and Proposition \ref{tilting theorem}.
\end{proof}

We also need the following observation for abelian categories.

\begin{lemma}\label{TnA}
Let $\AA$ be an abelian category.
\begin{enumerate}
\item Any object in $\Morph^{\sm}_{N-2}(\Proj\AA)$ is projective in $\Morph_{N-2}(\AA)$.
\item If $\PP$ is a subcategory of $\AA$ of projective generators, then $\Morph_{N-2}^{\sm}(\PP)$ is a subcategory of $\Morph_{N-2}(\AA)$ of projective generators.
\end{enumerate}
Assume that $\AA$ is an $Ab3$-category with a small full subcategory {$\CC$} of compact projective generators.
\begin{enumerate}
\item[(3)] $\Morph_{N-2}(\AA)$ is an $Ab3$-category with a small full subcategory $\Morph^{\sm}_{N-2}(\CC)$ of compact projective generators.
\item[(4)] We have an equivalence $\Morph_{N-2}(\AA)\simeq\Mod(\Morph^{\sm}_{N-2}(\CC))$ given by
$X\mapsto\Hom_{\Morph_{N-2}(\AA)}(-,X)|_{\Morph^{\sm}_{N-2}(\CC)}$.
\end{enumerate}
\end{lemma}

\begin{proof}
(1) 
By Lemma \ref{find tilting}(1), it suffices to prove that $\mu_i^{N-1}(C)$ is projective in 
$\Morph_{N-2}(\AA)$ for $C\in\Proj\AA$ and $1\le i\le N-1$. 
Indeed, let an epimorphism $Y\to X$ in $\Morph_{N-2}(\AA)$ be given. 
Then it induces an epimorphism 
$\Hom _{\AA} (C, Y^{N-i} ) \to \Hom _{\AA} (C, X^{N-i} ) $. 
Since $X^{N-i}=\opn{H}^{N-i}_{(i)}(X)$ and $Y^{N-i}=\opn{H}^{N-i}_{(i)}(Y)$, 
we get an epimorphism 
$\Hom _{\KKK _N (\AA )} (\mu_i^{N-1}(C), Y) \to \Hom _{\KKK _N (\AA )} (\mu_i^{N-1}(C), X) $ from Lemma \ref{homfrommu}(2). 
\par\noindent
(2) Let $X=(X^1\xrightarrow{\alpha^1}\cdots\xrightarrow{\alpha^{N-2}}X^{N-1})$ be any object in 
$\Morph_{N-2}(\AA)$. For each $1\le i\le N-1$, we take an epimorphism $P_i\to X^i$ with $P_i\in\PP$. Then we have an epimorphism $\coprod_{i=1}^{N-1}\mu^{N-1}_{N-i}(P_i)\to X$.
\par\noindent
(3) The assertion follows from (1), (2) and Lemma \ref{find tilting}(4).
\par\noindent
(4) This is immediate from (3) and Proposition \ref{Aus2}.
\end{proof}

Now we are ready to prove Theorem \ref{DND}.

\begin{proof}[of Theorem \ref{DND}]
By Proposition \ref{find tilting2} and Lemma \ref{TnA}, we have triangle equivalences $\DDD_N(\AA)\simeq\DDD(\Mod(\Morph_{N-2}^{\sm}(\CC)))\simeq\DDD(\Morph_{N-2}(\AA))$,
which restrict to the identity functor on $\Morph_{N-2}^{\sm}(\CC)$. 
\end{proof}

Next, to restrict the above equivalence to the subcategories of bounded complexes, we give the following preliminary result.

\begin{lemma}\label{detect cohomology}
Let $\AA$ be an abelian category and {$\CC$} a full subcategory of projective generators.
Then the following conditions are equivalent for $X\in\DDD_N(\AA)$.
\begin{enumerate}
\item $X$ belongs to $\DDD_N^{\bo}(\AA)$ (resp., $\DDD_N^-(\AA)$, $\DDD_N^+(\AA)$).
\item For every $0<r<N$, $\Hom_{\DDD_N(\AA)}(\mu^s_r(\CC),X)=0$ holds for all but finitely many (resp., sufficiently large, sufficiently small) $s\in\z$.
\item $\Hom_{\DDD_N(\AA)}(\Morph_{N-2}^{\sm}(\CC),\Sigma^iX)=0$ holds for all but finitely many (resp., sufficiently large, sufficiently small) $i\in\z$.
\end{enumerate}
\end{lemma}

\begin{proof} 
(1) and (2) are equivalent by Lemma \ref{homfrommu}(2).
\par\noindent
Since $\Sigma^2=\Theta^N$ holds and $\DDD_N^{\bo}(\AA)$ (resp., $\DDD_N^-(\AA)$, $\DDD_N^+(\AA)$) is closed under $\Sigma$, the condition (2) is equivalent to the following condition.
\begin{itemize}
\item For any $0<r<N$ and $0\le s<N$, $\Hom_{\DDD_N(\AA)}(\mu^s_r(\CC),\Sigma^iX)=0$ holds for all but finitely many (resp., sufficiently large, sufficiently small) $i\in\z$.
\end{itemize}
This is equivalent to the condition (3) since $\tri\{\mu^s_r(P)\mid P\in\CC,\ 0<r<N, 0\le s<N\}=\KKK_N^{\bo}(\CC)=\tri\Morph_{N-2}^{\sm}(\CC)$ holds 
by Lemmas \ref{Sigma of mu}(2) and \ref{find tilting}(3).
\end{proof}

Now we are able to prove the following result.

\begin{theorem}\label{DND2}
Let $\AA$ be an $Ab3$-category with a small full subcategory of compact
projective generators. 
Then the triangle equivalence in Theorem \ref{DND} restricts to those for $\natural=+,-,\bo$
\[
\DDD^{\natural}_{N}(\AA)\simeq \DDD^{\natural}(\Morph_{N-2}(\AA)).\]
\end{theorem}

\begin{proof}
This is immediate from Theorem \ref{DND} and Lemma \ref{detect cohomology}.
\end{proof}

In the case $\mcal{A} = \Mod {R}$ for a ring $R$, 
$\Morph_{N-2}(\AA)$ is nothing but  the category of modules over the upper triangular matrix ring  $\opn{T}_{N-1}(R)$ of size $N-1$ over $R$.  
Then we have the following precise description of homologies. 

\begin{proposition}\label{describe H by H}
Let $R$ be a ring. 
Then we have a triangle equivalence \[ G:\DDD_{N}(\Mod R)\simeq \DDD(\Mod \opn{T}_{N-1}(R)) \] 
which gives the following for $X\in\DDD_N(\Mod R)$ and $i\in\z$:
\begin{eqnarray*}
\opn{H}^{2i}(GX) &=& \left(\opn{H}_{(N-1)}^{iN+1}(X) \to \opn{H}_{(N-2)}^{iN+2}(X) \to \cdots \to \opn{H}_{(1)}^{iN+N-1}(X)\right), \\
\opn{H}^{2i+1}(GX) &=& \left(\opn{H}_{(1)}^{(i+1)N}(X) \to \opn{H}_{(2)}^{(i+1)N}(X) \to \cdots \to \opn{H}_{(N-1)}^{(i+1)N}(X)\right),
\end{eqnarray*}
where each morphism is a canonical one between homologies.
\end{proposition}

\begin{proof}
By Theorem \ref{DND}, we have a triangle equivalence $G:\DDD_{N}(\Mod R)\simeq \\ \DDD(\Mod \opn{T}_{N-1}(R))$ which is the identity on $\Morph^{\sm}_{N-2}(\pj R)$.
We shall show the equalities only for $i=0,1$ since for others it follow from $\Theta^N=\Sigma^2$. 
For $0<r<N$, we have
\begin{eqnarray*}
\Hom_{\Mod \opn{T}_{N-1}(R)}(\mu_{r}^{N-1}(R), \opn{H}^{0}(GX))&\simeq &
\Hom_{\KKK (\Mod \opn{T}_{N-1}(R)) }(\mu_{r}^{N-1}(R), GX)\\
\simeq \Hom_{\DDD(\Mod \opn{T}_{N-1}(R))}(\mu_{r}^{N-1}(R), GX)
&\simeq&\Hom_{\DDD_{N}(\Mod R)}(\mu_{r}^{N-1}(R), X)\simeq \opn{H}_{(r)}^{N-r}(X). 
\end{eqnarray*}
The first isomorphism is from Lemma \ref{TnA}(1), the second from 
$\mu_{r}^{N-1}(R) \in \KKK_N^{\mrm{p}}(\Proj R)$, {{and the}}
the third by $G$. The last is from Lemma \ref{homfrommu}(2).
Thus the morphism $\opn{H}_{(r+1)}^{N-r-1}(X) \to \opn{H}_{(r)}^{N-r}(X)$ is the canonical one since it is induced from the canonical morphism $\mu_{r}^{N-1}(R) \to \mu_{r+1}^{N-1}(R)$.
Similarly we have
\begin{eqnarray*}
&&\Hom_{\Mod \opn{T}_{N-1}(R)}(\mu_{r}^{N-1}(R), \opn{H}^{1}(GX))\simeq\Hom_{\DDD(\Mod \opn{T}_{N-1}(R))}(\Sigma^{-1}\mu_{r}^{N-1}(R), GX) \\
&\simeq&\Hom_{\DDD_{N}(\Mod R)}(\Sigma^{-1}\mu_{r}^{N-1}(R), X)\simeq \Hom_{\DDD_{N}(\Mod R)}(\mu_{N-r}^{N-r-1}(R), X)\simeq\opn{H}_{(N-r)}^0(X)
\end{eqnarray*} as desired.
\end{proof}

As an application, we have the following results for homotopy categories.

\begin{corollary}\label{KNhtp}
Let $\BB$ be an additive category with arbitrary coproducts.
If $\BB^{\co}$ is skeletally small and satisfies $\BB=\Sum(\BB^{\co})$,
then we have triangle equivalences
\[\KKK^{-}_{N}(\BB)\simeq \KKK^{-}(\Morph^{\sm}_{N-2}(\BB))\ \mbox{ and }\ 
\KKK^{\bo}_{N}(\BB)\simeq\KKK^{\bo}(\Morph^{\sm}_{N-2}(\BB)).\]
\end{corollary}

\begin{proof}
Let $\AA=\Mod\BB^{\co}$. Then $\AA$ (resp., $\Morph_{N-2}(\AA)$) is an $Ab3$-category with 
a subcategory $\BB$ (resp., $\Morph_{N-2}^{\sm}(\BB)$) of projective generators by Lemma \ref{TnA}(2). Thus we have triangle equivalences
\[\KKK^{-}_{N}(\BB)\simeq\DDD_{N}^{-}(\AA)\simeq\DDD^{-}(\Morph_{N-2}(\AA))
\simeq\KKK^{-}(\Morph^{\sm}_{N-2}(\BB)).\] 
where the first and the third equivalence by Remark \ref{rem:spcpx} and the second by Theorem \ref{DND2}.
Since these equivalences restrict to the identity functor on $\Morph^{\sm}_{N-2}(\BB)$, we have a triangle equivalence
\[\KKK^{\bo}_{N}(\BB)=\tri_{\KKK^{-}_{N}(\BB)}\Morph^{\sm}_{N-2}(\BB)\simeq
\tri_{\KKK^{-}(\Morph^{\sm}_{N-2}(\BB))}\Morph^{\sm}_{N-2}(\BB)=\KKK^{\bo}(\Morph^{\sm}_{N-2}(\BB))\]
by Lemma \ref{find tilting}(3).
\end{proof}

\begin{example}\label{GrMod1}
Let $R$ be a graded ring, and $\cat{GrMod}R$ the category of graded right $R$-modules.
Then $\cat{GrMod}R$ satisfies the condition of Theorem \ref{DND}. Hence we have
a triangle equivalence for $\natural=$nothing$,-,\bo$:
\[\DDD_{N}^{\natural}(\cat{GrMod}R)\simeq\DDD^{\natural}(\Morph_{N-2}(\cat{GrMod}R)).\]
\end{example}

Finally we study the bounded derived category of $N$-complexes in the case of coherent rings. We prepare the following easy observation.

\begin{lemma}\label{Gmu}
Let $G:\DDD_N(\AA)\to\DDD(\Morph_{N-2}(\AA))$ be the triangle equivalence given in Theorem \ref{DND}. 
For any $P\in\CC$ and $i, r \in\z$ with $0\le r<N$, we have
\[G(\mu_1^{iN+r}(P))=\left\{\begin{array}{ll}
\cdots\to0\to{\mu^{N-1}_{N-1}(P)}\to0\to\cdots&\mbox{if }\ r=0,\\
\cdots\to0\to{\mu^{N-1}_{N-r-1}(P)}\to{\mu^{N-1}_{N-r}(P)}\to0\to\cdots&\mbox{if }\ 0<r<N.
\end{array}\right.\] 
which is a complex concentrated in degree $2i-1$ if $r=0$, in $2i-1$ and $2i$ otherwise. 
\end{lemma}

\begin{proof}
Since $\Sigma^2=\Phi^N$, we have only to show them for the case $i=0$ by an induction on $r$. 
If $r=0$, then we have $G(P)=\Sigma\mu^{N-1}_{N-1}(P)$ since $P=\Sigma\mu^{N-1}_{N-1}(P)$.
Assume $0<r<N$. 
Then an exact sequence $0\to \mu^{N-1}_{N-r-1}(P)\to\mu^{N-1}_{N-r}(P)\to\mu^r_1(P)\to 0$ in $\CCC _N (\AA)$ induces 
a triangle $\mu^{N-1}_{N-r-1}(P)\to\mu^{N-1}_{N-r}(P)\to\mu^r_1(P)\to\Sigma\mu^{N-1}_{N-r-1}(P)$ in $\DDD _N (\AA)$ 
by Proposition \ref{prop:dertrig}(1).
Applying $G$, we have a triangle $\mu^{N-1}_{N-r-1}(P)\to\mu^{N-1}_{N-r}(P)\to G\mu^r_1(P)\to\Sigma\mu^{N-1}_{N-r-1}(P)$ in $\DDD_N(\AA)$.
\end{proof}

\begin{proposition}\label{DNMod}
Let $R$ be a ring.
\begin{enumerate}
\item We have triangle equivalences for $\natural=-, \bo, (-,\bo)$:
\[
\KKK^{\natural}_{N}(\pj R)\simeq \KKK^{\natural}(\pj \opn{T}_{N-1}(R)).
\]
\item If $R$ is right coherent, then we have triangle equivalences for $\natural=-,\bo$:
\[\DDD^{\natural}_{N}(\fmod R)\simeq\DDD^{\natural}(\fmod \opn{T}_{N-1}(R)).\]
\end{enumerate}
\end{proposition}

\begin{proof}
(1)  According to Theorem \ref{cor:subeqv}, we regard
$\KKK^-_{N}(\Pj R)$ (resp., $\KKK^{-}(\Pj \opn{T}_{N-1} (R)$)
as a full subcategory of $\DDD_N(\Mod R)$ (resp., $\DDD(\Mod {T}_{N-1} (R))$).   
We shall show that the triangle equivalence 
$G:\DDD_N(\Mod R)\simeq\DDD(\Mod \opn{T}_{N-1}(R))$ in Theorem \ref{DND} 
restricts to the desired equivalence. 
Indeed, $G$ induces a triangle equivalence
\[\begin{aligned}
\KKK^{\bo}_{N}(\pj R) & =\tri_{\DDD_N(\Mod R)}\Morph_{N-2}^{\sm}(\pj R)\simeq\tri_{\DDD(\Mod \opn{T}_{N-1}(R))}\pj\opn{T}_{N-1}(R) \\ &=\KKK^{\bo}(\pj \opn{T}_{N-1}(R)).
\end{aligned}\]
\par\noindent
To get the triangle equivalence for $\natural = {-}$,  
we shall show $GP\in  \KKK^-(\pj \opn{T}_{N-1}(R))$ for each $P \in \KKK^{-}_{N}(\pj R)$.
We may assume $P \in \CCC^{-}_{N}(\pj R)$ and $\tau _{\geq 1} P=0$. Set $P_{n}=\tau_{\geq -n}P$ for each $n>0$.
Then we have a term-wise split exact sequence $0 \to P_{n-1} \to P_{n} \to \Theta^{n} P^{-n} \to 0$ in $\CCC^{\bo}_{N}(\pj R)$, 
and a triangle in $\DDD_{N}(\Mod R)$
\[
P_{n-1} \to P_{n} \to \Theta^{n} P^{-n} \stackrel{\varphi _n} {\to }\Sigma P_{n-1}.
\]
Applying $G$, we have a triangle in $\DDD(\Mod\opn{T}_{N-1}(R))$
\[
GP_{n-1} \to GP_{n} \to G\Theta^{n} P^{-n} \stackrel{G \varphi _n} {\to } \Sigma GP_{n-1}.
\]
There exists a term-wise split exact sequence
\[
0\to Q_{n-1}\to Q_n\to G\Theta^nP^{-n}\to0
\]
in $\CCC^{\bo}(\pj \opn{T}_{N-1}(R))$ such that $GP_0\to GP_1\to GP_2\to\cdots$ is isomorphic to $Q_0\to Q_1\to Q_2\to\cdots$.
Then Lemma \ref{Gmu} gives a triangle $GP_{n-1}\to GP_n\to G\Theta^nP^{-n}\to\Sigma GP_{n-1}$ such that 
$G\Theta^{n}P^{-n}$ has only non-zero terms at degrees $2\lfloor n/N\rfloor$ and $2\lfloor n/N\rfloor-1$, where $\lfloor n/N\rfloor$ is 
the largest integer $m$ satisfying $m\le n/N$. 
Therefore  $\tau_{>2\lfloor n/N\rfloor}Q_{n-1}=\tau_{>2\lfloor n/N\rfloor}Q_{n}$ 
hence $\underset{\longrightarrow}{\lim}\ Q_n\in\KKK^{-}(\pj \opn{T}_{N-1}(R))$.
Since $P\simeq\underset{\longrightarrow}{\opn{hlim}}\ P_n$ in $\DDD_{N}(\Mod R)$ by Lemma \ref{prop:limithlimit}, 
$GP\simeq \underset{\longrightarrow}{\opn{hlim}}\ GP_n\simeq\underset{\longrightarrow}{\lim}\ Q_n$ in 
$\DDD(\Mod\opn{T}_{N-1}(R))$. Thus $GP\in\KKK^{-}(\pj \opn{T}_{N-1}(R))$ holds.
\par\noindent
By a similar argument, a quasi-inverse functor $G^{-1}:\DDD(\Mod\opn{T}_{N-1}(R))\simeq\DDD_N(\Mod R)$ induces a functor $\KKK^-(\pj\opn{T}_{N-1}(R))\simeq\KKK_N^-(\pj R)$.
Hence $G$ restricts to a triangle equivalence $\KKK_{N}^-(\pj R)\simeq \KKK^-(\pj\opn{T}_{N-1}(R))$.
By Lemma \ref{detect cohomology}, this restricts to 
a triangle equivalence $\KKK_{N}^{-,\bo}(\pj R)\simeq \KKK^{-,\bo}(\pj\opn{T}_{N-1}(R))$.
\par\noindent
(2) When $R$ is right coherent, $\opn{T}_{N-1}(R)$ is also right coherent.
In fact, let $A$ be $\opn{T}_{N-1}(R)$ and $e_i$ ($1\le i\le N-1$) the idempotent of $A$ whose $(i,i)$-entry is $1$ and others are zero. 
Let $0\to Z\to Y\to X$ be an exact sequence of $A$-modules such that $X$ and $Y$ are finitely presented.
Since $e_iAe_i=R$, we have an exact sequence $0\to  Ze_i\to  Ye_i\to  Xe_i$ of $R$-modules.
The $R$-modules $Xe_i$ and $Ye_i$ are finitely presented and
$R$ is coherent, hence so is the $R$-module $Ze_i$ for any $1\le i\le N-1$.
Therefore the $A$-module $Z$ is finitely generated.

We have the desired triangle equivalences
\begin{eqnarray*}
&\DDD^{-}_{N}(\fmod R) \simeq \KKK^{-}_{N}(\pj R)\simeq\KKK^{-}(\pj \opn{T}_{N-1}(R))
\simeq\DDD^{-}(\fmod \opn{T}_{N-1}(R)),& \\
&\DDD^{\bo}_{N}(\fmod R) \simeq \KKK^{-, \bo}_{N}(\pj R)\simeq\KKK^{-, \bo}(\pj \opn{T}_{N-1}(R)) \simeq\DDD^{\bo}(\fmod \opn{T}_{N-1}(R))&
\end{eqnarray*}
from (1) for the middles, Theorem \ref{cor:subeqv} for the others. 
\end{proof}



\end{document}